\documentclass{article}
\usepackage{amsmath,amsfonts,amssymb,amsthm,bbm,latexsym,mathrsfs}
\usepackage{graphicx,color,epsfig,fancyhdr,dsfont}
\usepackage{enumerate}
\usepackage{hyperref}
\usepackage{indentfirst}
\usepackage{amsmath,amscd}
\usepackage{caption}
\usepackage{graphicx, subfig}
\usepackage[all]{xy}

\newtheorem{definition}{Definition}
\newtheorem{lemma}[definition]{Lemma}

\newtheorem{coro}[definition]{Corollary}
\newtheorem{prop}[definition]{Proposition}

\newtheorem*{claim}{Claim}

\newenvironment{claimproof}[1][\proofname]
               {
                 \proof[#1]
                 
               }
               {
                 \endproof
               }

\newcommand{\Lie}{\operatorname{Lie}}

\newcommand{\BC}{\operatorname{BC}}

\newcommand{\Cert}{\operatorname{Cert}}
\newcommand{\Good}{\operatorname{Grp}}

\newcommand{\Zcl}{\operatorname{Zcl}}
\newcommand{\cl}{\operatorname{cl}}

\newtheorem{maintheorem}{Theorem}

\newcommand{\pCF}{p\operatorname{CF}}
\newcommand{\gG}{\mathfrak{g}}
\newcommand{\tp}{\operatorname{tp}}

\begin{document}
\title{On the uniformity of definably amenable radicals and dfg-fsg decompositions  in \(p\)-adically closed fields}

\author{Zhentao Zhang}
\date{}
\maketitle

\begin{abstract}
For definable families of groups over $p$-adically closed fields, we show that the definably amenable radicals and their dfg-fsg decompositions can be defined uniformly as the prime $p$ and the family parameters vary.
\end{abstract}

\section{Introduction}

Two structural decomposition results have been established for groups definable in $p$-adically closed fields. In \cite{YZ}, Yao and Zhang  proved the existence of definably amenable radicals and obtained the corresponding definably amenable--semisimple decompositions.
In \cite{PYZ}, Pillay, Yao and Zhang established dfg--fsg decompositions for definably amenable groups. In this paper, we study uniform versions of these results as the prime and the parameters of a fixed definable family vary.

A definable group is \emph{definably amenable} if there is a finitely additive probability measure on its definable subsets, that is invariant under left translation. In an NIP theory, a definable group is definably amenable if and only
if it admits a global strongly $f$-generic type
\cite[Corollary~5.10]{HP2011}.
Two important special cases are dfg groups and fsg groups. A definable group has \emph{dfg} (definable $f$-generics) if some global complete type on it has every left translate definable over a common small model. It has \emph{fsg} (finitely satisfiable generics) if some global complete type has every left translate finitely satisfiable in a common small model. Over a $p$-adically closed field, fsg is equivalent to definable compactness \cite{JohnsonFsg}. Dfg groups, after passage to a finite-index subgroup, admit definable homomorphisms with finite kernel and finite-index image into split solvable algebraic groups \cite[Theorem~5.6]{PYZ}.

Every definably amenable group $G$ over a $p$-adically closed field
admits a definable exact sequence
$$1\longrightarrow H\longrightarrow G\longrightarrow C\longrightarrow1$$
where $H$ has dfg and $C$ has fsg \cite[Theorem~8.8]{PYZ}. Such a sequence is called a
\emph{dfg--fsg decomposition}, and its kernel $H$ is a
\emph{normal dfg component} of $G$.

More generally, a definable dfg subgroup $H\leq G$ is called a
\emph{dfg component} of $G$ if it has maximum dimension among all
definable dfg subgroups of $G$. This is equivalent to maximality
under almost inclusion among such subgroups, and also to definable
compactness of the homogeneous space $G/H$
\cite[Theorems~6.4 and~8.18]{PYZ}. Here $H_1$ is almost contained
in $H_2$ if $H_1\cap H_2$ is of finite index in $H_1$. When $G$ is definably
amenable, its dfg components are pairwise commensurable
\cite[Corollary~8.9]{PYZ}, and the normal ones are precisely the
kernels occurring in dfg--fsg decompositions.

Every group $G$ definable in a $p$-adically closed field $M$ has a greatest
definable normal definably amenable subgroup, called its
\emph{definably amenable radical} and denoted by $D(G)$
\cite[Proposition~3.1.17]{YZ}. The radical $D(G)$ is uniquely
determined and definable over the same parameters as $G$.
It yields a \emph{definably amenable-semisimple decomposition}
of $G$, given by a definable exact sequence
$$1\longrightarrow D(G)\longrightarrow G
 \xrightarrow{\pi} S\longrightarrow1$$
where $S$ has a definable subgroup of finite index which is definably isomorphic to a finite-index subgroup of
$\mathbf S(M)$ for some connected semisimple algebraic group $\mathbf S$ over $M$ \cite[Corollary~3.1.10]{YZ}.
In particular, $S$ is trivial when $G$ is definably amenable.

The preceding results hold within each fixed $p$-adically closed
theory, where definable amenability is also definable in families
\cite[Proposition~2.4.11]{YZ}. We study the uniformity of these
constructions as the prime varies. We now state our main results more precisely.

Let $L=\{0,1,+,-,\cdot\}$ and $\pCF$ be the theory of $\mathbb {Q}_p$ in $L$ for each prime $p$. Fix formulas
\[
 \gamma(x;w),\qquad \mu(x,y,z;w),\qquad \iota(x,y;w)
\]
where $\gamma$ specifies the underlying set, while $\mu$ and $\iota$
specify multiplication and inversion, respectively.
Let $\Good(w)$ express that these data define a group.
For $M=M_p\models\pCF$, a tuple $a$ in $M$ is called
\emph{admissible} if $M\models\Good(a)$. We denote the resulting group with underlying set $\gamma(M;a)$ by $G_{M,a}$. 

\begin{maintheorem}\label{main-A}
There is an $L$-formula $\delta_D(x;w)$ such that
$\delta_D(M_p;a)=D(G_{M_p,a})$ for every prime $p$, $M_p\models \pCF$ and admissible $a$ in $M_p$.
\end{maintheorem}

Note that $G$ is definably amenable iff $G=D(G)$.
As a corollary, for given $\gamma(x;w)$, $\mu(x,y,z;w)$ and $\iota(x,y;w)$, there is an $L$-formula  $\Good^*(w):=\Good(w)\wedge \forall x(\gamma(x;w)\leftrightarrow \delta_D(x;w))$ expressing that these data define a definably amenable group. For convenience, we call a tuple $a$ in $M=M_p\models \pCF$ \emph{strongly admissible} if $M\models \Good^*(a)$. 

Moreover, we can uniformly realize the corresponding
definably amenable--semisimple decompositions as definable exact sequences.

\begin{maintheorem}
\label{main-AS}
There are  $L$-formulas such that,
for every prime $p$, every $M_p\models\pCF$, and every admissible
$a$ from $M_p$, there is an additional parameter tuple $b$ 
making these formulas define a definably amenable--semisimple decomposition
 $$1\longrightarrow D(G_{M_p,a})\longrightarrow G_{M_p,a}
 \xrightarrow{\pi_{a,b}}S_{M_p,a,b}\longrightarrow1.$$

Moreover, the additional parameters $b$ are used only to specify
a definable Lie atlas on $G_{M_p,a}$. In particular, if the given
family consists of definable subgroups of $\mathrm{GL}_m$ for some fixed $m$, then the exact sequences can be defined uniformly over $a$ alone.
\end{maintheorem}

For dfg-fsg decomposition for definably amenable groups, we have

\begin{maintheorem}\label{main-H}
There is an $L$-formula $\delta_H(x;w)$ such that
$\delta_H(M_p;a)$ is a normal dfg component of $G_{M_p,a}$ for every prime $p$, $M_p\models \pCF$ and strongly admissible $a$ in $M_p$.
\end{maintheorem}

The corresponding dfg-fsg decompositions can be also uniformly realized.

\begin{maintheorem}
\label{main-HC}
There are  $L$-formulas such that,
for every prime $p$, every $M_p\models\pCF$, and every strongly admissible
$a$ from $M_p$, there is an additional parameter tuple $b$ 
making these formulas define a dfg-fsg decomposition
 $$1\longrightarrow K_{M_p,a}\longrightarrow G_{M_p,a}
 \xrightarrow{\rho_{a,b}}C_{M_p,a,b}\longrightarrow1.$$
Moreover, if the given family consists of subgroups of
$\mathrm{GL}_m$ for some fixed $m$, the dfg-fsg decompositions
can be defined uniformly over $a$ alone, without additional
parameters.
\end{maintheorem}

Section~\ref{sec:preliminaries} recalls definable compactness and
gives a uniform description of the Lie group structures on the
groups in the family. In Section~\ref{sec:radical}, we use the
adjoint action on the semisimple quotient of the Lie algebra to
prove Theorems~\ref{main-A} and~\ref{main-AS}.
Section~\ref{sec:dfg-preparations} develops growth chains in
pseudo-local fields and gives uniform codes for Zariski closures.
In Section~\ref{sec:uniform-dfg-fsg}, we use these results to
choose normal dfg components uniformly and to construct the
quotient maps required for Theorems~\ref{main-H} and~\ref{main-HC}.

The language $L$ and the formulas defining the family are fixed
throughout. For each prime $p$, we write $M=M_p\models\pCF$.
Unless otherwise stated, definability allows parameters from the field under consideration.

We use boldface letters for algebraic groups. For an algebraic
group $\mathbf G$ over $M$, we write $\mathbf G(M)$ for its group
of $M$-rational points and $\operatorname{Lie}(\mathbf G)$ for
the $M$-rational Zariski tangent space at the identity, equipped
with its natural Lie bracket. For $X\subseteq M^n$, we write $\Zcl(X)$
for its Zariski closure in affine $n$-space over an algebraic
closure of $M$. All Lie algebras are finite-dimensional
over the base field, and all Lie algebra homomorphisms and
derivations are linear over that field.

\section{Uniform Lie group structures}\label{sec:preliminaries}

We recall definable compactness and study the uniformity of
definable Lie group structures in the given family.

\subsection{Valuation and definable compactness}

By \cite{CDLM}, there is a parameter-free $L$-formula $o(z)$ defining the valuation
ring $\mathcal {O}_M$ in every $M=M_p\models\pCF$, uniformly in $p$. The maximal ideal of $\mathcal{O}_M$ is therefore uniformly definable by
\[
 \mathfrak{m} (z):=(z=0)\vee
 \bigl(o(z)\wedge\exists y\,(zy=1\wedge\neg o(y))\bigr).
\]
We write
\[
 v:M^\times\longrightarrow\Gamma_M
\]
for the valuation, with $v(0)=\infty$. Absolute values take values
in the associated multiplicative ordered value group, and
\[
 \|x\|=\max_i|x_i|
\]
for tuples $x$ from $M$. Unless otherwise stated, affine spaces are equipped with the valuation topology.
Comparisons between absolute values, and the corresponding
neighborhood conditions, are expressible in $L$ using $o$.

A definable set $X\subseteq M^r$ is \emph{bounded} if
$cX\subseteq\mathcal O_M^r$ for some $c\in M^\times$.
A definable topological space is \emph{definably compact} if every
definable downward directed family of nonempty closed subsets has
nonempty intersection.

For a definable subset of affine space with its induced topology,
definable compactness is equivalent to being closed and bounded
\cite[Fact~2.4]{PYZ}. Continuous definable images, closed definable
subspaces, finite products and finite unions preserve definable
compactness \cite{JY-1}. In particular, a finite-index extension of
a definably compact definable group is definably compact.

For matrix groups, we use the closed embedding
\[
 \mathrm{GL}_d(M)\longrightarrow M^{d^2+1}:
 g\longmapsto\bigl(g,\det(g)^{-1}\bigr).
\]
Thus boundedness in $\mathrm{GL}_d(M)$ includes boundedness of the
matrix entries and of the inverse determinant.

Every definable group carries a canonical definable group topology
\cite{HP,JY-1}, which we describe using definable Lie atlases in the
next subsection. This topology may differ from the topology induced
by a given presentation of the group in affine space.
Definable compactness of groups is always understood with respect
to their group topology.

\subsection{Definable Lie atlases}\label{R:sec:atlas}

We recall the weak Lie atlases of
\cite[Definitions~5.1--5.5]{ALH}.
The definitions and first-order descriptions below also apply to
any nontrivially valued field of characteristic zero in which
$o$ defines the valuation ring.

A definable map
$r:U\longrightarrow M^d$ with $U\subseteq M^l$
is \emph{topologically \'etale} if it is a local homeomorphism,
where $U$ carries the induced topology.
A finite weak atlas on a definable group $G$ consists of finitely many definable
injections
\[
 \phi_j:U_j\longrightarrow G
\]
whose images cover $G$, together with definable topologically
\'etale maps
\[
 r_j:U_j\longrightarrow M^d.
\]
We require the overlap domains
\[
 U_{ij}:=\phi_i^{-1}\bigl(\phi_j(U_j)\bigr)
\]
to be open in $U_i$.

Using local inverses of the $r_j$, we express transition maps and
group operations in field coordinates.
The atlas is a \emph{weak strictly differentiable $T_3$ Lie atlas}
if these maps are strictly differentiable and satisfy the $T_3$
condition of \cite[Definition~3.11]{ALH}.
The charts induce the topology on $G$.

An \emph{atlas template} $\mathcal T$ is a finite collection of
$L$-formulas intended to define the chart domains $U_j$ and the
graphs of $r_j$ and $\phi_j$.
These formulas have parameter variables $(w,v)$, where $w$
consists of the original group-family parameters and $v$ consists
of additional atlas parameters.
Substituting tuples $(a,b)$ gives candidate atlas data for
$G_{M,a}$, with the length of $b$ fixed by $\mathcal T$.

\begin{lemma}
\label{R:lem:certificate}
For every atlas template $\mathcal T$, there is an $L$-formula
$\Cert_{\mathcal T}(w,v)$ such that, for every prime $p$, every
$M=M_p\models\pCF$ and all tuples $a,b$ from $M$ corresponding
to $w,v$, the condition
\[
 M\models\Cert_{\mathcal T}(a,b)
\]
holds iff $a$ is admissible and the data specified by
$\mathcal T$ with parameters $(a,b)$ form a weak strictly
differentiable $T_3$ Lie atlas on $G_{M,a}$.

The same formula expresses the corresponding conditions in
every nontrivially valued field of characteristic zero in which
$o$ defines the valuation ring.
\end{lemma}

\begin{proof}
Functionality, injectivity, coverage, relative openness and the
local-homeomorphism conditions are first order for the specified
graphs.

For a map $f$ between open subsets of affine spaces, strict
differentiability at $z$ with derivative $A$ is expressed by the
condition that, for every $\epsilon\ne0$, there is $\eta\ne0$
such that
$$
 \|f(u)-f(v)-A(u-v)\|
 \leq |\epsilon|\,\|u-v\|
$$
whenever $\|u-z\|<|\eta|$ and $\|v-z\|<|\eta|$.

The $T_3$ condition is also first order. Locally, it requires
degree-three Taylor expansions with a fourth-order remainder
bound uniform in the center of expansion. The coefficient of
degree $r$ must itself admit a polynomial approximation of degree
$3-r$ at the given point.
All degrees and dimensions are fixed, so the coefficients and
bounds are finite tuples of field elements.

Taylor coefficients are uniquely determined: the difference
between two candidate expansions can be restricted to lines,
where uniqueness follows by comparing its lowest nonzero term
with the remainder bound.
Consequently, the coefficient functions have definable graphs,
which can be substituted into the preceding conditions.

Local inverse branches of the $r_j$ are described by reversing
their graphs on sufficiently small neighborhoods.
Applying these descriptions to the transition maps,
multiplication and inversion gives $\Cert_{\mathcal T}$.
\end{proof}

We briefly recall the associated Lie algebra.
Let $e=1_G$ be the identity of $G$, and put
\[
 \gG=\Lie(G):=T_eG.
\]
A chart containing $e$ identifies $\mathfrak g$ with $M^d$.
For $g\in G$, let $c_g:G\to G$ be the conjugation map
$c_g(x)=gxg^{-1}$. The adjoint action is defined by
$\operatorname{Ad}(g):=\mathrm{d}(c_g)_e$.
The chain rule makes
\[
 \operatorname{Ad}:G\longrightarrow
 \operatorname{Aut}_M(\mathfrak g)
\]
a definable homomorphism. Its derivative at $e$ gives the bracket
\[
 [X,Y]=\mathrm{d}(\operatorname{Ad})_e(X)(Y).
\]
This is a Lie bracket
\cite[Propositions~6.16 and~6.19]{ALH}.

For a fixed atlas template $\mathcal T$, these constructions are
uniformly definable from $(a,b)$.
Indeed, choose the first chart containing the identity $e=1_{G_{M,a}}$ of $G=G_{M,a}$.
The preimage of $e$ in that chart is uniquely determined, and
derivatives and Taylor coefficients are defined by the conditions
used in the preceding proof.
Thus, the structure constants of $\mathfrak g$ and the matrices
of $\operatorname{Ad}$ require no further parameters.

Definable subgroups are closed and carry the topology induced
by the ambient group \cite[Fact~2.12]{JY-3}.
Definable homomorphisms are differentiable, and surjective ones
are submersions and hence open. Consequently, the group topology
on any definable realization of a group quotient agrees with
the quotient topology.

We also record a standard fact for linear groups.
Suppose that $G\leq\mathrm{GL}_d(M)$, and let $\mathbf {G}$
be its Zariski closure. 
By \cite{Mac}, $\dim G=\dim\mathbf G(M)$, so $G$ has
nonempty interior in $\mathbf G(M)$.
Then $G$ is open in $\mathbf G(M)$ and
 $\operatorname{Lie}(G)=\operatorname{Lie}(\mathbf G)$.

\subsection{Uniformity across primes}

We now show that finitely many atlas templates suffice for the
given family, uniformly over all primes.

\begin{lemma}[Uniform existence of Lie atlases]
\label{R:lem:finitecover}
There are finitely many atlas templates
$\mathcal T_1,\ldots,\mathcal T_s$ such that, for every prime $p$,
every $M=M_p\models\pCF$ and every admissible tuple $a$ from $M$,
$$M\models
 \bigvee_{i=1}^s\exists v_i
 \Cert_{\mathcal T_i}(a,v_i).$$
\end{lemma}

\begin{proof}
Enumerate all atlas templates as
$\mathcal T_1,\mathcal T_2,\ldots$.
For each fixed prime $p$, every definable group admits a
definable Lie atlas \cite{HP}.
Thus every admissible tuple in every model of $\pCF$ admits
an atlas represented by one of these templates.
By compactness, a finite initial segment suffices for all
models of this fixed theory.

Suppose that no finite initial segment works uniformly in $p$.
The fixed-prime conclusion allows us to choose $p_j>j$,
models $M_{p_j}$ and admissible tuples $a_j$ such that for $i\leq j$,
\[
 M_{p_j}\models
 \neg\exists v_i \Cert_{\mathcal T_i}(a_j,v_i).
\]

Let $N=\prod_j M_{p_j}/\mathcal U$ be a nonprincipal ultraproduct and $a=(a_j)/\mathcal U$.  
Then $N\models\Good(a)$, but for every $i$,
\[
 N\models
 \neg\exists v_i\Cert_{\mathcal T_i}(a,v_i).
\]

Equip $N$ with the valuation ring defined by $o$.
Since $p_j\to\infty$, this is a nontrivially valued Henselian
field of equicharacteristic zero.
Its pure valued-field theory is $1$-h-minimal
\cite[Corollary~6.2.6]{CHRK}.
In a sufficiently saturated elementary extension $N'$,
the group $G_{N',a}$ therefore admits a finite weak strictly
differentiable $T_3$ Lie atlas
\cite[Proposition~6.4]{ALH}.

Since the valuation ring is defined by $o$, the atlas data
are given by $L$-formulas with finitely many field parameters.
They are consequently represented by some template
$\mathcal T_i$ and a tuple $b$ from $N'$, so $ N'\models\Cert_{\mathcal T_i}(a,b)$.
By elementarity, we have
$N\models\exists v_i\Cert_{\mathcal T_i}(a,v_i)$, a contradiction.
\end{proof}

Since the list of atlas templates is finite, we may use a single
parameter tuple $b$ of fixed length. Its first coordinate records
the chosen template, and the remaining coordinates give its
parameters, with unused coordinates set to zero.
The associated Lie algebra and adjoint action are then uniformly
definable over $(a,b)$.

\section{Uniform definably amenable radicals}\label{sec:radical}

We prove the uniform definability of the radical and construct
the associated definably amenable--semisimple decompositions.

Recall that definable amenability is preserved under definable
homomorphic images, extensions, finite-index subgroups and finite
extensions. Moreover, finite, solvable and definably compact
definable groups are definably amenable
\cite{HPP,JohnsonFsg,PYZ}.

\subsection{Bounded conjugacy}

Let $A\leq\mathrm{GL}_d(M)$ be a definable subgroup, and let
\[
 \BC(A)=
 \{u\in A:\{huh^{-1}:h\in A\}\text{ is bounded}\}.
\]
Since the determinant is constant on each conjugacy class,
$\BC(A)$ is defined by
\[
 u\in A\ \wedge\
 \exists c\ne0\ \forall h\in A\
 \bigwedge_{1\leq i,j\leq d}
 o\bigl(c(huh^{-1})_{ij}\bigr).
\]
Thus the subgroups $\BC(A)$ are uniformly $L$-definable in
definable families of matrix groups, using only the parameters
defining $A$.
Boundedness of conjugacy classes is preserved by a fixed
change of basis.

We recall the following facts about semisimple groups.
Let $\mathbf S$ be a connected $M$-simple algebraic group.
If $\mathbf S$ is $M$-isotropic, then the subgroup
$\mathbf S(M)^+$ generated by the $M$-points of the unipotent
radicals of its $M$-defined parabolic subgroups is definable
and of finite index in $\mathbf S(M)$
\cite[Corollary~3.15(i)]{PYZ}, and hence Zariski dense.
It is abstractly simple \cite[Section~3.2]{PYZ}
and is not definably amenable \cite[Lemma~4.3]{PYZ}.
In this case, every definable open subgroup of $\mathbf S(M)$
which is not definably compact contains $\mathbf S(M)^+$
\cite[Corollary~3.15(ii)]{PYZ}.
If $\mathbf S$ is $M$-anisotropic, then $\mathbf S(M)$ is
definably compact \cite[Corollary~2.17]{PYZ}.

For a finite-dimensional semisimple $M$-Lie algebra $\mathfrak s$,
let $\operatorname{Aut}(\mathfrak s)$ denote its algebraic
automorphism group and $\operatorname{Int}(\mathfrak s)$ its
connected algebraic group of inner automorphisms.
Since every derivation of $\mathfrak s$ is inner
\cite[Chapter~I, Theorem~4.22]{Milne},
\[
 \operatorname{Lie}(\operatorname{Aut}(\mathfrak s))
 =\operatorname{Der}(\mathfrak s)
 =\operatorname{ad}(\mathfrak s)
 =\operatorname{Lie}(\operatorname{Int}(\mathfrak s)).
\]
As $\operatorname{Int}(\mathfrak s)$ is a connected algebraic
subgroup of $\operatorname{Aut}(\mathfrak s)$, it follows that
$\operatorname{Int}(\mathfrak s)=\operatorname{Aut}(\mathfrak s)^\circ$
where $\circ$ denotes the Zariski identity component.

This group is semisimple, since its Lie algebra is isomorphic
to $\mathfrak s$ \cite[Chapter~II, Proposition~4.1]{Milne}.
Moreover, $\operatorname{ad}$ identifies its faithful action on
$\mathfrak s$ with its adjoint representation, so its center
is trivial.

\begin{lemma}\label{R:lem:semisimplebc}
Let $\mathfrak s$ be a finite-dimensional semisimple $M$-Lie
algebra. Let
$A\leq\operatorname{Aut}(\mathfrak s)(M)$ be a definable subgroup satisfying
 $\operatorname{Lie}(A)=\operatorname{ad}(\mathfrak s)$.
Then $\BC(A)$ is a definably compact normal subgroup of $A$
containing every definable normal definably amenable subgroup
of $A$.
\end{lemma}

\begin{proof}
Let
$\mathbf S=\operatorname{Int}(\mathfrak s)$ and $
 A_0=A\cap\mathbf S(M)$.
Since $\operatorname{Aut}(\mathfrak s)/\mathbf S$ is finite,
$A_0$ is normal and of finite index in $A$.
Its Lie algebra is $\operatorname{Lie}(\mathbf S)$, so $A_0$
is open in $\mathbf S(M)$.

The semisimple group $\mathbf S$ is an almost direct product
of connected almost $M$-simple algebraic subgroups
\cite[Fact~2.12]{PYZ}.
Since $\mathbf S$ is centerless, this is a direct product
\[
 \mathbf S=\mathbf S_1\times\cdots\times\mathbf S_r,
\]
where each $\mathbf S_i$ is connected and $M$-simple.

Let $A_i$ be the projection of $A_0$ to $\mathbf S_i(M)$.
Then $A_i$ are definable open subgroups of
$\mathbf S_i(M)$.

Let $I$ be the set of indices $i$ for which $A_i$ is not
definably compact, and let $J$ be the set of indices $i$
for which $A_i$ is definably compact.
For $i\in I$, the group $\mathbf S_i$ is $M$-isotropic and 
$P_i:=\mathbf S_i(M)^+\leq A_i$ is of finite index.
We regard $P_i$ as a subgroup of $\mathbf S(M)$ via the inclusion of the $i$-th factor. Since $A_0$ is open in $\mathbf S(M)$,
$A_0\cap P_i$ is open in $P_i$ and hence nontrivial. It is normal
in $P_i$, since every element of $P_i\leq A_i$ lifts to an element of $A_0$ inducing the same conjugation on the $i$-th factor.
By simplicity of $P_i$, $A_0\cap P_i=P_i$, so $P_i\leq A_0$.

Let
\[
 \mathbf S_I=\prod_{i\in I}\mathbf S_i \qquad\text{ and }
 \qquad
 \mathbf S_J=\prod_{i\in J}\mathbf S_i.
\]
Conjugation by $A$ permutes the factors $\mathbf S_i$
and the corresponding subgroups $A_i$, preserving definable
compactness. Hence $\mathbf S_I$ and $\mathbf S_J$ are
normalized by $A$.

\begin{claim}
 $\BC(A)=C_A(\mathbf S_I(M))$ in
$\operatorname{Aut}(\mathfrak s)$.
\end{claim}

\begin{claimproof}
Let $u\in\BC(A)$ and $i\in I$.
For an $M$-defined one-parameter unipotent subgroup
 $v:\mathbb G_a\longrightarrow\mathbf S_i$,
we have $v(M)\subseteq P_i\subseteq A$.
In the representation on $\mathfrak s$, the entries of
$v(t)uv(t)^{-1}$
are polynomials in $t$ and are bounded on $M$.
A nonconstant polynomial is unbounded on $M$, since its leading
term dominates for arguments of sufficiently negative valuation.
Thus, these entries are constant.
Evaluating at $t=0$ shows that $u$ centralizes $v(M)$.
These one-parameter subgroups generate the Zariski-dense
subgroup $P_i$, and hence $u$ centralizes $\mathbf S_I(M)$.

Conversely, $C_A(\mathbf S_I(M))$ is normal in $A$, since
$\mathbf S_I(M)$ is normalized by $A$.
The direct product decomposition
$\mathbf S=\mathbf S_I\times\mathbf S_J$
and the triviality of $Z(\mathbf S)$ imply that
$ C_A(\mathbf S_I(M))\cap A_0
 =\mathbf S_J(M)\cap A_0$.
Since $A_0$ is closed in $\mathbf S(M)$,  $C_A(\mathbf S_I(M))\cap A_0$
is closed in $\mathbf S_J(M)$.
Then it is a closed
subgroup of a definably compact group, so it is definably compact.
Since $A_0$ is of finite index in $A$, $C_A(\mathbf S_I(M))\cap A_0$ is of
finite index in $C_A(\mathbf S_I(M))$.
Thus, $C_A(\mathbf S_I(M))$ is definably compact.
By normality, every $A$-conjugacy class of one of its elements
is contained in this group and is therefore bounded.
Consequently, $C_A(\mathbf S_I(M))\subseteq\BC(A)$.
\end{claimproof}

Let $L\lhd A$ be a definably amenable subgroup.
For $i\in I$, let $L_i$ be the projection of $L\cap A_0$
to $\mathbf S_i(M)$.
Then $L_i\lhd A_i$ is definably amenable.
Since $P_i$ has finite index in $A_i$, the subgroup
$L_i\cap P_i$ is definably amenable.
Since $P_i$ is simple and is not definably amenable, $L_i\cap P_i$ is
trivial. 
Since $P_i$ is Zariski dense in $\mathbf S_i$ and the
centralizing condition is Zariski closed, we have $L_i=1$. Hence, $L\cap A_0\leq \mathbf S_J(M)$.

For $u\in L$ and $s\in\prod_{i\in I}P_i$, the commutator
$usu^{-1}s^{-1}$ belongs to both $L\cap A_0\leq \mathbf S_J(M)$ and
$\mathbf S_I(M)$, and is therefore trivial.
Since $P_i$ is Zariski dense in $\mathbf S_i$ and the
centralizing condition is Zariski closed, we have 
$L\leq C_A(\mathbf S_I(M))=\BC(A)$.
\end{proof}

\subsection{The radical and its linear quotient}

Let $G$ be a group definable in $M$.
We apply Lemma~\ref{R:lem:semisimplebc} to the action of $G$ on the
semisimple quotient of its Lie algebra.
This gives a description of $D(G)$ and a definable linear
realization of $G/D(G)$.

Let
\[
 \mathfrak g=\operatorname{Lie}(G),\qquad
 \mathfrak r=\operatorname{rad}(\mathfrak g),\qquad
 \mathfrak s=\mathfrak g/\mathfrak r.
\]
Every automorphism of $\mathfrak g$ preserves its solvable
radical $\mathfrak r$. Hence, the adjoint action induces
a definable homomorphism
\[
 \alpha:G\longrightarrow\operatorname{Aut}(\mathfrak s)(M).
\]
Let $A=\alpha(G)$. Since $\alpha:G\to A$ is a submersion and $\mathrm{d}\alpha_{1_G}(X)=\operatorname{ad}(X+\mathfrak r)$
for every $X\in\mathfrak g$, we have
$\operatorname{Lie}(A)=\mathrm{d}\alpha_{1_G}(\mathfrak g)
 =\operatorname{ad}(\mathfrak s)$.

\begin{prop}\label{R:prop:localradical}
 $D(G)=\alpha^{-1}(\BC(\alpha(G)))$.
\end{prop}

\begin{proof}
Let $L\lhd G$ be any definable normal definably amenable subgroup.
Its image $\alpha(L)$ is a definable normal definably amenable
subgroup of $A:=\alpha(G)$. 

Then Lemma~\ref{R:lem:semisimplebc} gives
$\alpha(L)\leq\BC(A)$.
Thus, $L\leq\alpha^{-1}(\BC(A))$.

Now it suffices to show that $\alpha^{-1}(\BC(\alpha(G)))$ is definably amenable. Since $\BC(\alpha(G))$ is definably amenable, it suffices to show that $\ker\alpha$ is definably amenable.

The subgroup $\ker\operatorname{Ad}$ is commutative-by-finite
\cite[Lemma~2.11 and its proof]{JY-3}, and therefore definably
amenable. Let $K=\operatorname{Ad}(G)$, and let $\mathbf K$ be its
Zariski closure in $\operatorname{Aut}(\mathfrak g)$.
By the linear-group observation in the preceding section, we have
$ \operatorname{Lie}(\mathbf K)
 =\operatorname{Lie}(K)
 =\operatorname{ad}(\mathfrak g)$.
Every automorphism of $\mathfrak g$ preserves $\mathfrak r$,
so the action of $\mathbf K$ on $\mathfrak s$ is an algebraic
homomorphism. Let $\mathbf H$ be its kernel.
Since $\mathfrak s$ has trivial center, we have
$$\operatorname{Lie}(\mathbf H)=\{\operatorname{ad}(X):
   X\in\mathfrak g \text{ with } [X,\mathfrak g]\subseteq\mathfrak r\}
 =\operatorname{ad}(\mathfrak r).$$
This Lie algebra is solvable.
It follows that $\mathbf H^\circ$ is solvable
\cite[Chapter~II, Corollary~4.4]{Milne}.
Thus, $K\cap\mathbf H(M)$ is solvable-by-finite and hence, definably amenable.
Since $\ker\alpha$ is an extension of $K\cap\mathbf H(M)$ by $\ker\operatorname{Ad}$, it is definably amenable.
\end{proof}

We next realize the quotient by $D(G)$ as a linear definable
group, using the same adjoint action.

\begin{coro}\label{R:cor:quotientlocal}
There is a definable surjective homomorphism
$\pi:G\to S$ onto a linear definable group $S$ with
$\ker\pi=D(G)$.
The group $S$ has no nontrivial definable normal definably
amenable subgroup.
\end{coro}

\begin{proof}
Let 
$\mathfrak c=
 \operatorname{ad}^{-1}\bigl(\operatorname{Lie}(\BC(A))\bigr)$ and
 $\mathfrak q=\mathfrak s/\mathfrak c$.

With the notation of the proof of
Lemma~\ref{R:lem:semisimplebc}, we have
$\BC(A)=C_A(\mathbf S_I(M))$ and 
 $\operatorname{Lie}(\BC(A))=\operatorname{Lie}(\mathbf S_J)$.
Thus, $\mathfrak c$ is the $A$-invariant sum of the simple
ideals corresponding to $J$.
The natural projection
$\mathfrak s\to\mathfrak q=\mathfrak s/\mathfrak c$
restricts to an $A$-equivariant isomorphism from the
complementary sum corresponding to $I$ onto $\mathfrak q$.
By the connectedness of $\mathbf S_I$ and the injectivity
of $\operatorname{ad}$, an element of $A$ acts trivially
on $\mathfrak q$ iff it centralizes
$\mathbf S_I(M)$.
Hence, the induced action
\[
 \pi:G\longrightarrow\operatorname{Aut}(\mathfrak q)(M)
\]
has kernel $\ker\pi
 =\alpha^{-1}(C_A(\mathbf S_I(M)))
 =\alpha^{-1}(\BC(A))
 =D(G)$.
Then let $S:=\pi(G)$, and it is a linear definable group.

Let $L\lhd S$ be definably amenable normal subgroup.
Then $\pi^{-1}(L)$ is a definable normal definably
amenable subgroup of $G$ containing $D(G)$.
Then $\pi^{-1}(L)=D(G)$, and hence $L$ is trivial.
\end{proof}

\subsection{Uniform formulas and quotient maps}
\label{sec:uniformradical}\label{R:sec:descent}

We now express the preceding constructions uniformly in the
given family. For a fixed atlas template, all constructions
below use only the group parameters and the atlas parameters.

\begin{lemma}\label{R:lem:linearalgebra}
Let $(\mathfrak g_t)_{t\in T}$ be a family of Lie algebras
such that the underlying linear subspaces
$\mathfrak g_t\subseteq M^d$ and their Lie brackets are
uniformly definable, where $d$ is fixed.
Then the derived algebras $[\mathfrak g_t,\mathfrak g_t]$
and the solvable radicals $\operatorname{rad}(\mathfrak g_t)$
are uniformly definable.

Moreover, if a uniformly definable family of group actions
on $\mathfrak g_t$ preserves a uniformly definable family
of ideals $\mathfrak i_t$, then the induced actions on
$\mathfrak g_t/\mathfrak i_t$ can be represented in matrix
coordinates uniformly and definably.
All these constructions use only the parameters defining
the given families and actions.
\end{lemma}

\begin{proof}
Since $d$ is fixed, Gaussian elimination to reduced
row-echelon form uses only a bounded number of field
operations and zero tests. Hence, bases of $\mathfrak g_t$
can be chosen uniformly on $t$. The derived algebra
$[\mathfrak g_t,\mathfrak g_t]$ is the span of the brackets
of these basis vectors, so it is uniformly definable.

For $X\in\mathfrak g_t$, membership in
$\operatorname{rad}(\mathfrak g_t)$ is equivalent to the
existence of at most $d$ vectors whose span is a solvable
ideal containing $X$.
The ideal condition and containment of $X$ are first order.
Since the dimension is at most $d$, solvability is equivalent
to the vanishing of the $d$-th derived algebra.
Iterating the bracket-span construction $d$ times expresses
this condition uniformly.
Existential quantification over the spanning vectors
therefore uniformly defines
$\operatorname{rad}(\mathfrak g_t)$.

For the final assertion, row reduction in the chosen
coordinates gives uniformly definable bases of
$\mathfrak g_t/\mathfrak i_t$.
The entries of the induced matrices are the uniquely
determined coordinates of the images of these basis
vectors under the given actions.
They are consequently uniformly definable using only
the given parameters.
\end{proof}

We also need a uniform description of the Lie algebra of
a matrix group. Recall that $\mathfrak m(z)$ denotes the
fixed $L$-formula defining the maximal ideal of the
valuation ring.

\begin{lemma}\label{R:lem:tangent}
Let $G\leq\mathrm{GL}_d(M)$ and $X\in M^{d\times d}$. Then 
$X\in\operatorname{Lie}(G)$ iff
$$\forall\epsilon\ne0\
  \exists t\ne0\ \exists h\in G
  \left(\mathfrak m(t/\epsilon)\ \wedge\
 \bigwedge_{1\leq i,j\leq d}
 \mathfrak m\left(
 \frac{(h-I-tX)_{ij}}{\epsilon t}
 \right)
 \right).$$
The Lie bracket is $[X,Y]=XY-YX$.
Thus, the Lie algebras of matrix groups are uniformly
$L$-definable from the parameters defining the groups.
\end{lemma}

\begin{proof}
The equivalence is the usual characterization of the tangent
space at $1_G=I$ by limits of $(h-I)/t$, as $t\to0$ through
nonzero values and $h\to I$ in $G$.
The displayed formula expresses this condition using
$\mathfrak m$, and is therefore uniform in $p$ and in the
parameters defining $G$.
The bracket is the usual matrix commutator.
\end{proof}

Now we prove Theorem \ref{main-A}.

\begin{proof}[Proof of Theorem~\ref{main-A}]
Fix an atlas template $\mathcal T$.
For every $M=M_p\models\pCF$ and all tuples $a,b$ satisfying
$M\models\Cert_{\mathcal T}(a,b)$, let
$\mathfrak g_{a,b}$ be the coordinate realization of
$\operatorname{Lie}(G_{M,a})$ supplied by the atlas.
By Section~\ref{R:sec:atlas}, this Lie algebra and its
adjoint action are uniformly definable over $(a,b)$.
Lemma~\ref{R:lem:linearalgebra} then uniformly defines
\[
 \mathfrak r_{a,b}
 =\operatorname{rad}(\mathfrak g_{a,b})
 \qquad \text{and} \qquad
 \mathfrak s_{a,b}
 =\mathfrak g_{a,b}/\mathfrak r_{a,b}
\]
together with the induced homomorphism
\[
 \alpha_{a,b}:G_{M,a}\longrightarrow
 \operatorname{Aut}(\mathfrak s_{a,b})(M).
\]

Since bounded conjugacy is uniformly $L$-definable,
there is an $L$-formula $P_{\mathcal T}(x;w,v)$ such that
$P_{\mathcal T}(x;a,b)$ defines
\[
 \alpha_{a,b}^{-1}
 \bigl(\BC(\alpha_{a,b}(G_{M,a}))\bigr)
\]
in every valid instance.
By Proposition~\ref{R:prop:localradical}, this subgroup
is $D(G_{M,a})$.

Let $\mathcal T_1,\ldots,\mathcal T_s$ be the finite list
provided by Lemma~\ref{R:lem:finitecover}.
Define
$$\delta_D(x;w):=\gamma(x;w)\ \wedge \bigvee_{i=1}^s\exists v_i
 \left(
 \Cert_{\mathcal T_i}(w,v_i)
 \wedge P_{\mathcal T_i}(x;w,v_i)
 \right).$$
Every admissible $a$ admits a valid atlas from this list,
and every valid atlas gives the same subgroup $D(G_{M,a})$.
Hence,
\[
 \delta_D(M;a)=D(G_{M,a})
\]
for every prime $p$, every $M=M_p\models\pCF$ and every
admissible $a$.
\end{proof}

As a corollary, definable amenability is definable uniformly on $p$.

\begin{coro}\label{R:cor:amenability}
For the given family, the set of admissible parameters $a$
such that $G_{M,a}$ is definably amenable is uniformly
$L$-definable for all primes $p$ and all
$M=M_p\models\pCF$.
\end{coro}

\begin{proof}
Since $G_{M,a}$ is definably amenable iff
$D(G_{M,a})=G_{M,a}$, Theorem~\ref{main-A} shows that
the required set is defined by
\[
 \Good^*(w):=\Good(w)\wedge
\forall x\bigl(\gamma(x;w)\leftrightarrow\delta_D(x;w)\bigr).
\]
\end{proof}

Now we prove Theorem~\ref{main-AS}.

\begin{proof}[Proof of Theorem~\ref{main-AS}]
Fix an atlas template $\mathcal T$ from
Lemma~\ref{R:lem:finitecover}. Let $M=M_p\models\pCF$ and tuples $a,b$ satisfy
$M\models\Cert_{\mathcal T}(a,b)$.
Use the Lie algebras $\mathfrak g_{a,b}$,
$\mathfrak r_{a,b}$, $\mathfrak s_{a,b}$ and the
homomorphism $\alpha_{a,b}$ constructed in the proof
of Theorem~\ref{main-A}.

Let $A_{a,b}=\alpha_{a,b}(G_{M,a})$ and define
\[
 \mathfrak c_{a,b}
 =\operatorname{ad}^{-1}
   \bigl(\operatorname{Lie}(\BC(A_{a,b}))\bigr)
 \qquad \text{and} \qquad
 \mathfrak q_{a,b}
 =\mathfrak s_{a,b}/\mathfrak c_{a,b}
\]
where $\operatorname{ad}$ is the adjoint isomorphism
$\mathfrak s_{a,b}\to\operatorname{ad}(\mathfrak s_{a,b})$.
The proof of Corollary~\ref{R:cor:quotientlocal} shows that
$\mathfrak c_{a,b}$ is an $A_{a,b}$-invariant ideal and
that the induced action of $G_{M,a}$ on $\mathfrak q_{a,b}$
has kernel $D(G_{M,a})$.

By the uniform definition of bounded conjugacy and
Lemma~\ref{R:lem:tangent},
$\operatorname{Lie}(\BC(A_{a,b}))$ is uniformly definable
over $(a,b)$.
Lemma~\ref{R:lem:linearalgebra} therefore gives coordinates
on $\mathfrak q_{a,b}$ and the matrices of the induced
action uniformly over $(a,b)$.
Since
\[
 \dim\mathfrak q_{a,b}\leq\dim G_{M,a}\leq n:= \text{the length of }x,
\]
we extend these matrices by identity blocks to obtain
homomorphisms into $\mathrm{GL}_n(M)$.

Encoding the finitely many atlas templates and their
parameters by a tuple $b$ of fixed length, as described
after Lemma~\ref{R:lem:finitecover}, gives uniformly
definable homomorphisms
\[
 \pi_{a,b}:G_{M,a}\longrightarrow\mathrm{GL}_n(M).
\]
Let $S_{M,a,b}=\pi_{a,b}(G_{M,a})$.
We obtain definable exact sequences
\[
 1\longrightarrow D(G_{M,a})
 \longrightarrow G_{M,a}
 \xrightarrow{\pi_{a,b}}S_{M,a,b}
 \longrightarrow1.
\]
The targets and their group operations are uniformly
definable as matrix groups.
By Corollary~\ref{R:cor:quotientlocal}, these are the required
definably amenable-semisimple decompositions.

The additional parameters $b$ specify only the atlas,
since all subsequent constructions are uniformly
definable from its Lie algebra and adjoint action.
In particular, if the given family consists of subgroups
of $\mathrm{GL}_m$ for a fixed $m$, the matrix coordinates
make all constructions uniformly definable over $a$ alone.
\end{proof}

\section{Subgroup chains and definable quotients}
\label{sec:dfg-preparations}

This section establishes the technical results needed for the uniform construction of dfg components and their quotients. We first describe subgroup chains by first-order
conditions. We then construct one-dimensional subgroups and definable quotients in suitable pseudo-local fields. Finally, we prove a uniform degree bound for linear families.

\subsection{Growth chains}
\label{D:subsec:chains}

We work in $M=M_p\models\pCF$. Let $G$ be a definable group. Recall that
definable subgroups are closed and carry the induced group topology \cite[Fact~2.12]{JY-3}.

Write $\Gamma$ for the value group, interpreted using the valuation ring formula $o$.
Following \cite[Proposition~2.28]{JY-1}, $G$ admits a \emph{$\Gamma$-exhaustion} $\{E_\lambda:\lambda\in\Gamma\}$, namely, a
 neighborhood basis at $1_G$ satisfying that
\begin{itemize}
 \item Each $E_\lambda$ is definably compact and open,
 and $E_\lambda^{-1}=E_\lambda$.
 \item $E_{\lambda}\subseteq E_{\lambda'}$ when $\lambda\leq\lambda'$.
 \item $\displaystyle\bigcup_{\lambda\in\Gamma}E_\lambda=G$.
\end{itemize}
For a fixed Lie atlas template, the closed-ball construction uses fixed $L$-formulas over the group and atlas parameters.
Thus, the corresponding $\Gamma$-exhaustions are uniformly definable in families.

A subset of $G$ is \emph{bounded} if it is contained
in some $E_\lambda$.
Products and translates of bounded sets are bounded.
For definable subgroups $H\leq V\leq G$, 
we call the quotient $V/H$ \emph{bounded} if
$V\subseteq E_\lambda H$ for some $\lambda\in\Gamma$. If the quotient $V/H$ is definable, then it is bounded
iff it is definably compact with the quotient topology. Recall that when $H\lhd V$ and $V/H$ is a definable group, the quotient topology agrees with the canonical
group topology. A complete type is \emph{unbounded} if every definable
set belonging to it is unbounded.

\begin{definition}
A \emph{growth chain} in $G$ is a chain of definable subgroups
\[
 \{1_G\}=H_0\lhd H_1\lhd\cdots\lhd H_r=H\leq G
\]
such that $\dim H_i=i$ for $0\leq i\leq r$ and
$H_{i+1}/H_i$ is unbounded for $0\leq i<r$.
A growth chain is \emph{complete} if $G/H$ is bounded.
\end{definition}

Recall that the orbit space of a definable action of a dfg
group on a definable set is definable
\cite[Theorem~4.1]{AGJ}.

We show that every term $H_i$ of a growth chain has dfg. The case $H_0=\{1_G\}$ is immediate.
Assume that $H_i$ has dfg. Then $H_{i+1}/H_i$ is definable. It is one-dimensional and unbounded, so it has dfg \cite[Lemma~6.10]{JY-1}.
Then $H_{i+1}$ is an extension of dfg groups, so it has dfg.

If the chain is complete, $G/H$ is definable and bounded,
hence definably compact. Thus, $H$ is a dfg component of $G$
\cite[Theorem~8.18]{PYZ}.
If, in addition, $H\lhd G$, then $G/H$ has fsg.

The following normal-core construction will be used to
make the endpoint of a growth chain normal in a larger group.
Suppose a growth chain ends at $H$, and let $H\lhd U\lhd V\leq G$ be
definable subgroups such that
 $U/H$ is bounded.
Then $H$ has dfg and $U/H$ has fsg. The normal core
\[
 \operatorname{Core}_V(H)
 =\bigcap_{v\in V}vHv^{-1}
\]
has finite index in $H$ \cite[Lemma~8.2]{PYZ}.
Replacing each $H_i$ by $H_i\cap\operatorname{Core}_V(H)$
gives a growth chain of the same length whose endpoint
$\operatorname{Core}_V(H)$ is normal in $V$.
Moreover, $U/\operatorname{Core}_V(H)$ remains bounded,
and every subgroup of $H$ normal in $V$ is contained
in $\operatorname{Core}_V(H)$.

For fixed subgroup and exhaustion formulas, the growth-chain
conditions and the stated dimension and boundedness
conclusions are first order.
These formulations will be used in the ultraproduct arguments below.

\subsection{Growth chains in pseudo-local fields}
\label{D:subsec:auxiliary}

To establish uniform definability across primes,
we study the required constructions in pseudo-local fields. We let
\[
 N=\prod_j M_{p_j}/\mathcal U
\]
where $\mathcal U$ is a nonprincipal ultrafilter on $\mathbb N$ and $p_j\rightarrow \infty$.
Then $N$ is a pseudo-local field of residue characteristic
zero, with pseudofinite residue field and value group
a $\mathbb Z$-group.
The first-order conditions established in
Subsection~\ref{D:subsec:chains} hold in $N$.
For definable subsets of affine space, dimension agrees
with the dimension of the Zariski closure \cite{WY}.
In particular, zero-dimensional definable sets are finite.

For a definable complete type $q$ over $N$ and an elementary extension $N'\succeq N$, we write $q|N'$ for its canonical definable extension to $N'$.

For each $N'\succeq N$, let $k_{N'}$ and $\Gamma_{N'}$
denote its residue field and value group, respectively. Let $F_{N'}$ be the henselization of $N'(t^{\mathbb Q})$ with respect to the $t$-adic valuation.
We equip $F_{N'}$ with the composition of this valuation and the given valuation on $N'$.
Its residue field is $k_{N'}$ and its value group is
$\mathbb Q\mathbin{\mathrm{lex}}\Gamma_{N'}:=\mathbb Q\oplus\Gamma_{N'}$,
ordered lexicographically with the $\mathbb Q$-coordinate first.
Since $\lambda\mapsto(0,\lambda)$ is elementary in Presburger language and the residue-field embedding is the identity, the Ax-Kochen-Ershov principle gives
$N'\preceq F_{N'}$.

Every element $a\in F_{N'}$ bounded over $N'$ (i.e., $\tp(a/N')$ is bounded) has a unique
standard part $\operatorname{st}(a)\in N'$, given by its residue under the coarse $t$-adic valuation. It satisfies
\[
 v\bigl(a-\operatorname{st}(a)\bigr)>\Gamma_{N'}.
\]
Moreover, $F_{N'}$ is algebraic over $N'(t)$, so every finite tuple in $F_{N'}$ has transcendence degree at most
one over $N'$.

We also need definability of types in the ring language $L$.
We regard $F_{N'}$ as a structure in the auxiliary
three-sorted language $L_t$ of \cite{Pas}, using the
coarse $t$-adic valuation $v_t$, residue field $N'$,
value group $\mathbb Q$, and angular component
\[
 \operatorname{ac}_t(x)
 =\operatorname{res}_t\bigl(x t^{-v_t(x)}\bigr)
\]
for $x\neq0$, where $\operatorname{res}_t$ denotes the coarse residue map.

\begin{lemma}\label{D:lem:definabletypes}
For every $N'\succeq N$ and every finite tuple $a$ from
$F_{N'}$, the type $\tp(a/N')$ is definable, in $L$.
\end{lemma}

\begin{proof}
Fix an $L$-formula $\varphi(x;y)$.
By applying \cite[Theorem~4.1]{Pas} to the $L_t$-expansion
of $F_{N'}$, we reduce $\varphi(a;b)$ to conditions on
$v_t(P(a;b))$ and $\operatorname{ac}_t(P(a;b))$
for finitely many polynomials $P(x;y)$,
uniformly for all tuples $b$ from $N'$.

The relevant monomials in $a$ span a finite-dimensional
$N'$-vector space.
Successive cancellation of leading terms gives a basis
$e_1,\ldots,e_r$ satisfying
\[
 v_t(e_1)<\cdots<v_t(e_r)
 \qquad \text{and} \qquad\operatorname{ac}_t(e_i)=1.
\]
Then, for each relevant polynomial $P(x;y)$, there are
polynomials $c_{P,i}(y)\in N'[y]$ such that
\[
 P(a;b)=\sum_{i=1}^r c_{P,i}(b)e_i
\]
for every tuple $b$ from $N'$. Since the $e_i$ form a basis over $N'$, the expression
vanishes exactly when all $c_{P,i}(b)$ vanish.
Otherwise, if $i$ is the first index with $c_{P,i}(b)\ne0$, then
\[
 v_t(P(a;b))=v_t(e_i)
 \qquad \text{and} \qquad
 \operatorname{ac}_t(P(a;b))=c_{P,i}(b).
\]
On each of these finitely many cases, the value-group
conditions have fixed truth values, while the residue-field
conditions become $L$-formulas over $N'$.
Thus, $\varphi(a;y)$ is equivalent to an
$L(N')$-formula in $y$.
\end{proof}

As in the above proof, we shall use
\cite[Theorem~4.1]{Pas} in the auxiliary language $L_t$
to analyze $L$-definable sets.
Throughout, definability refers to $L$ unless otherwise stated.

Until the end of the proof of Lemma~\ref{D:lem:growth}, we let $H\lhd V\leq G$ be definable groups over $N$, with $Q=V/H$ unbounded, and let $\pi:V\to Q$ be the quotient map. We fix an $N$-definable $\Gamma$-exhaustion $(E_\lambda)_\lambda$ for $G$ in a monster model and let $\hat  {E}_\lambda=\pi(E_\lambda\cap V)$.
We equip $Q$ with the quotient topology.
The sets $\hat E_\lambda$ are closed and open,
form a neighborhood basis at $1_Q$, and exhaust $Q$. A subset of $Q$ is bounded if it is contained in some $\hat E_\lambda$.

\begin{lemma}\label{D:lem:escape}
There are an algebraic curve $C$ over $N$ and an
$N$-definable complete type $p$ concentrating on $C\cap V$
such that both $p$ and $q:=\pi(p)$ are unbounded.
Moreover, for every $N'\succeq N$, the types $p|N'$ and $q|N'$ are the unique unbounded complete extensions of $p$ and $q$, respectively, and $p|N'$ can be realized in $V(F_{N'})$.
\end{lemma}

\begin{proof}
We let
\[
 p_0^N(x)
 =
 \{x\ne0\}\cup \{ v(x)>\lambda:\lambda\in \Gamma_N\}
 \cup
 \{\exists z\,(z^k=x):k\geq1\}.
\]
By \cite[Propositions~2.2.2 and~2.3.2]{AcostaHils}, $p_0^N$ is complete. It is easy to see that the distinguished $t\in F_N$ realizes $p_0^N$ and so does each $t^{1/m}$ with $m\geq1$.
By Lemma~\ref{D:lem:definabletypes}, $p_0^N$ is definable, and its canonical definable extension to $N'\succeq N$
is $p_0^{N'}$.

By the unboundedness of $Q$, pick $a\in V(F_N)$
such that $\pi(a)\notin \hat E_{v(t)}$.
Since $v(t)>\Gamma_N$, both $a$ and $\pi(a)$ are
unbounded over $N$.
Moreover, $F_N$ is algebraic over $N(t)$, so
$\operatorname{trdeg}(a/N)=1$.
Thus, the Zariski closure $C$ of $a$ over $N$ is an algebraic curve.
Let $p=\tp(a/N)$ and $ q=\pi(p)$.
By Lemma~\ref{D:lem:definabletypes}, both types are $N$-definable. Obviously, they are unbounded.

Applying the local root construction of
\cite[Lemmas~3.5 and~3.7]{Pas} in an $N$-definable chart containing $a$, we obtain, for some $m\geq 1$, an $N$-definable
continuous map $f:D\to C\cap V$ with $f(t^{1/m})=a$, where $D$
is clopen in a sufficiently small punctured neighborhood of $0$.
Since $t^{1/m}\models p_0^N$, we have $p=f(p_0^N)$.

Let $N'\succeq N$, and $a'$ realize an unbounded
complete extension of $p$ over $N'$.
Then there is $u\models p^N_0$ such that $f(u)=a'$.
If $u$ were not infinitesimal over $N'$, it would lie in an $N'$-definable closed annulus of $D$ and then  $a'$ would be
bounded over $N'$, a contradiction.
Hence, $u\models p_0^{N'}$, and therefore
$a'\models p|N'$.

Finally, since $p=f(p_0^N)$ and $f$ is $N$-definable, we have $p|N'=f(p_0^{N'})$.
The type $p_0^{N'}$ is realized by $t^{1/m}\in F_{N'}$, so $p|N'$ is realized by
$f(t^{1/m})\in V(F_{N'})$. The uniqueness assertion for $q=\pi(p)$ follows immediately.
\end{proof}

We now adapt the bounded-gap and $\mu$-stabilizer arguments
of \cite[Sections~4--6]{JY-1} to $N$, using the treatment of interpretable groups in \cite[Section~5]{JY-2}.

For subsets $I,A\subseteq Q$, we write
\[
 I\diamond A
 :=\{x\in I:xA\cap I=\varnothing\}.
\]

\begin{lemma}\label{D:lem:boundedgaps}
Suppose there is an $N$-definable subgroup $W\lhd V$
containing $H$ such that $W/H$ is bounded and $V/W$
is commutative.
Let $I\subseteq G$ be $N$-definable such that
$\hat I:=\pi(I\cap V)$ is unbounded in $Q$.
Then for every
$\lambda\in\Gamma_N$, there is $\lambda'>\lambda$ such that
 $\hat I\diamond
 (\hat E_{\lambda'}\setminus\hat E_\lambda)$
is bounded.
\end{lemma}

\begin{proof}
For the fixed defining formulas, the hypotheses are first order and hence hold in $\mathcal U$-almost all
factors $M_{p_j}$ by {\L}o\'s's theorem.
Since the conclusion is also first order, it suffices
to prove it in these factors.

In each such factor, use the same notation for the
corresponding data and apply \cite[Lemma~5.6]{JY-2}
to the commutative quotient $V/W$.
Let $\theta:V/H\to V/W$ be the natural map.
Its kernel $W/H$ is bounded, so bounded sets have
bounded preimages.
For the index $\gamma>\lambda$ obtained in $V/W$,
choose $\lambda'>\gamma$ such that
\[
 \theta^{-1}(\theta(\hat E_\gamma))
 =\hat E_\gamma(W/H)\subseteq\hat E_{\lambda'}.
\]
Thus a difference whose image lies in
$\theta(\hat E_\gamma)\setminus\theta(\hat E_\lambda)$
lies in $\hat E_{\lambda'}\setminus\hat E_\lambda$.
The required gap set is therefore contained in the
preimage of the bounded gap set in $V/W$, and is bounded.
\end{proof}

We recall the formal stabilizer construction from
\cite[Section~2]{PS} and \cite[Section~5]{JY-1}.
Definable subsets of $Q$ and complete or partial types on $Q$ are understood in $L^{\mathrm{eq}}$. Of course, the above lemma can also be formulated in $L^{\mathrm{eq}}$.
Let $\mu_N$ be the partial type of infinitesimals at $1_Q$, and write $\operatorname{Stab}(\Sigma)$ for the stabilizer of a partial type $\Sigma$ under left translation.

For an $N$-definable complete type $q$ on $Q$, the partial type $\mu_N\cdot q$ is definable \cite[Claim~2.15]{PS}.
There are $N$-definable subgroups $K_\varphi\leq Q$ such that,
for every $N'\succeq N$,
\[
 \operatorname{Stab}\bigl(\mu_{N'}\cdot(q|N')\bigr)
 =\bigcap_\varphi K_\varphi(N').
\]

We also use standard parts on groups and quotients.
On the closed bounded chart pieces used to construct
$(E_\lambda)_\lambda$, coordinatewise residue for the coarse
$t$-adic valuation gives a standard part in $N'$.
Hence, every point of $V(F_{N'})$ bounded over $N'$ has
a unique standard part in $V(N')$.
A bounded point of $Q(F_{N'})$ has a bounded representative
in $V(F_{N'})$, whose standard part defines its standard
part in $Q(N')$ via $\pi$.
For $Y\subseteq F_{N'}^m$, we let $\operatorname{st}(Y)$
be the set of standard parts of the points of $Y$ which are bounded over $N'$.

\begin{lemma}\label{D:lem:stdimension}
Let $N'\succeq N$ and let $Y\subseteq F_{N'}^m$
be definable.
Then $\operatorname{st}(Y)$ is $N'$-definable and $\dim(\operatorname{st}(Y))\leq\dim(Y)$.
\end{lemma}

\begin{proof}
In $L_t$, we have 
$$ \operatorname{st}(Y)
 =\operatorname{res}_t
 \bigl(\{a\in Y:v_t(a_i)\geq0\text{ for all }i\}\bigr).$$
By \cite[Theorem~4.1]{Pas}, it is definable
in the pure residue field $N'$.

Let $\dim Y=d$ and $Z$ be the Zariski closure of $Y$. Then by \cite[Theorems~B(2),~C(3) and Fact~1.17(2)]{WY}, $\dim Z=d$.

For each set $J$ of $d+1$ coordinate indices, the projection of $Z$ onto these coordinates satisfies a nonzero polynomial relation
 $P_J\in F_{N'}[X_j:j\in J]$.
Divide $P_J$ by a coefficient of least $v_t$-value.
Then all its coefficients have nonnegative $v_t$-value,
and one coefficient is $1$. Applying $\operatorname{res}_t$ coefficientwise, we have 
a nonzero polynomial
$\operatorname{res}_t(P_J)\in N'[X_j:j\in J]$.
For every bounded $a\in Y$, $\operatorname{res}_t$ gives
\[
\operatorname{res}_t(P_J)\bigl(\operatorname{st}(a)_J\bigr)=0.
\]

Let $W$ be the algebraic set over $N'$ defined by all these equations $\operatorname{res}_t(P_J)=0$.
Then $\operatorname{st}(Y)\subseteq W(N')$.
Since no $d+1$ coordinates can be algebraically
independent in $W$,
$\dim W\leq d$, and hence
$\dim(\operatorname{st}(Y)) \leq\dim W \leq d
 =\dim Y$.
\end{proof}

\begin{lemma}\label{D:lem:growth}
Suppose there is an $N$-definable subgroup $W\lhd V$
containing $H$ such that $W/H$ is bounded and $V/W$
is commutative.
Then there is an $N$-definable subgroup $U\leq V$ containing $H$ such that
$U/H$ is unbounded and 
$$\dim U\leq\dim H+1.$$
If $H$ occurs in a growth chain, then
$\dim U=\dim H+1$.
\end{lemma}

\begin{proof}
Let $C$, $p$ and $q=\pi(p)$ be as in
Lemma~\ref{D:lem:escape}.
Let $N'\succeq N$ be $|N|^+$-saturated.
Choose $a\in V(F_{N'})$ realizing $p|N'$, and put
$b=\pi(a)$.
Then $b\models q|N'$ is unbounded over $N'$.

For each $N$-definable subset $I$ belonging
to $p$, put $\hat I=\pi(I\cap V)$ and let
$$X_I:=\{h\in  Q(N'): \hat I(F_{N'})\cap h\hat E_\lambda(F_{N'})b
 \neq\emptyset \text{ for every }\lambda\in\Gamma_{N'}\}.$$
Since $q$ is definable over $N$, $X_I$ is definable over $N$.
The argument of \cite[Lemma~5.9]{JY-2}, together with Lemma~\ref{D:lem:escape}, gives
\[
\operatorname{Stab}\bigl(\mu_{N'}\cdot(q|N')\bigr)
 =\bigcap_{I\in p}X_I.
\]

\begin{claim}
For every $I\in p$, the set $X_I$ is unbounded.
\end{claim}

\begin{claimproof}
Fix $I\in p$ and $\lambda\in\Gamma_{N'}$.
Applying Lemma~\ref{D:lem:boundedgaps} to $I^{-1}$ (note that $V/W$ is commutative),
there is $\lambda'>\lambda$ such that
$ (\hat I)^{-1}\diamond
 (\hat E_{\lambda'}\setminus\hat E_\lambda)$
is bounded.
Since $b^{-1}$ is unbounded over $N'$, it lies outside $ (\hat I)^{-1}\diamond
 (\hat E_{\lambda'}\setminus\hat E_\lambda)$.
So there is
$ c\in\hat E_{\lambda'}(F_{N'})
 \setminus\hat E_\lambda(F_{N'})$
such that $cb\in\hat I(F_{N'})$.

Since the annulus is closed, we have
$ \operatorname{st}(c)\in
 \hat E_{\lambda'}(N')\setminus\hat E_\lambda(N')$.
By definition of $X_I$, $\operatorname{st}(c)\in X_I$.
Since $\lambda\in\Gamma_{N'}$ is arbitrary,
$X_I$ is unbounded.
\end{claimproof}

It is easy to see that
$X_{I\cap I'}\subseteq X_I\cap X_{I'}$.
So every finite intersection of these $X_I$ is unbounded. Since the family $\{X_I\}_I$ has cardinality at most $|N|$, saturation shows that $\operatorname{Stab}(\mu_{N'}\cdot(q|N'))=\bigcap_{I\in p}X_I$ is unbounded.

It remains to show the dimension.
Let
 $I_C=C\cap V$ and $Y=I_Ca^{-1}H$. It is easy to check that $\dim Y\leq \dim H+1$.
For each chart $\phi_j:U_j\to G$, we let
\[
 Y_j=\{u\in U_j(F_{N'}):\phi_j(u)\in Y(F_{N'})\}.
\]
and let $Z_j$ be the closure of
$\operatorname{st}(Y_j)$ in the valuation topology.

Let $g\in V(N')$ satisfy $\pi(g)\in X_{I_C}$.
By definition, for every $\lambda\in\Gamma_{N'}$,
\[
 Y(F_{N'})\cap gE_\lambda(F_{N'})\neq\emptyset.
\]
This implies that
$\phi_j^{-1}(g)$ lies in $Z_j$ for every chart containing $g$.
Applying Lemma~\ref{D:lem:stdimension} in these charts, we have
\[
 \dim\pi^{-1}(X_{I_C})
 \leq \max_j\dim(Z_j)
 =\max_j\dim\operatorname{st}(Y_j)
 \leq\dim Y
 \leq 1+\dim H.
\]

Recall that 
$\operatorname{Stab}(\mu_{N'}\cdot(q|N'))=\bigcap_\varphi K_\varphi(N')$, where each
$K_\varphi\leq Q$ is definable over $N$.
Then $\bigcap_\varphi K_\varphi(N')=\bigcap_I X_I$. Since $N'$ is $|N|^+$-saturated, there is $\varphi$ such that $K:=K_\varphi(N')\subseteq X_{I_C}$. 
On the other direction, there is $I$ such that
$X_I\subseteq K$. So $K$ is unbounded.

Let $U=\pi^{-1}(K)$. Then $U$ is an $N$-definable subgroup of $V$ containing $H$.
The quotient $U/H$ is unbounded and
$\dim U\leq 1+\dim H$.
Hence, this $U$ is what we want.

Finally, suppose that $H$ occurs in a growth chain. For the fixed defining formulas, the chain conditions, the unboundedness of $U/H$, and $\dim U\leq\dim H+1$ hold in $\mathcal U$-almost all factors $M_{p_j}$ by {\L}o\'s's theorem. For convenience, denote the corresponding groups in the factors by $H'$ and $U'$. Then $H'$ has dfg.
If $\dim U'=\dim H'$, then $U'/H'$ is finite
by \cite[Corollary~4.4]{AGJ}, contradicting the
unboundedness of $U'/H'$.
Thus, $\dim U'=\dim H'+1$ holds in $\mathcal U$-almost all factors. Hence, by {\L}o\'s's theorem, $\dim U=\dim H+1$ in $N$.
\end{proof}

\begin{lemma}\label{D:lem:extension}
Under the hypotheses of Lemma~\ref{D:lem:growth},
suppose in addition that a growth chain ends at $H$.
Then this chain extends to a growth chain ending at
$K\lhd V$ such that $V/K$ is bounded.

If $V\lhd G$ and $H\lhd G$, the endpoint $K$ can also
be chosen normal in $G$.
\end{lemma}

\begin{proof}
Let $U$ be as in Lemma~\ref{D:lem:growth}. 
To normalize this extension, we let
\[
 U_0=\{u\in U: uwu^{-1}w^{-1}\in H\text{ for every }w\in W\}.
\]
In almost all factors, the corresponding quotients
$U'/H'$ and $W'/H'$ have dfg and fsg, respectively.
By \cite[Proposition~7.4]{PYZ}, the conjugation action
of $U'/H'$ on $W'/H'$ has finite image, so
$U_0'$ is of finite index in $U'$.
The result is first-order, so we have that
 $\dim U_0=\dim U$ and $U_0/H$ is unbounded over $N$.

Let $\widetilde H=\operatorname{Core}_V(U_0)$. Then $\widetilde H\lhd V$. Applying the statements in Subsection~4.1 to
$U'_0\lhd W'U'_0\lhd V'$ in $\mathcal U$-almost all factors, the corresponding quotients $U_0'/\widetilde {H}'$ are finite and $\dim\widetilde {H}'=\dim H'+1$. Hence, $U_0/\widetilde H$ is bounded and $\dim \widetilde H=\dim H+1$. It is easy to check that $\widetilde H$ extends the original growth chain. 

If $V/\widetilde H$ is bounded, we can let $K=\widetilde H$. If $V/\widetilde H$ is unbounded, we repeat with $(\widetilde H,W\widetilde H)$ in place of $(H,W)$. 
Since each step increases the dimension of the endpoint, iteration terminates at a normal endpoint $K$ in $V$ with bounded quotient.

If $V\lhd G$ and $H\lhd G$, we let $\widetilde K=\operatorname{Core}_G(K)$. In $\mathcal U$-almost all factors, the $G'$-conjugates
of $K'$ are normal dfg subgroups of $V'$ with fsg quotients, and hence are commensurable
\cite[Lemma~8.1]{PYZ}.
By \cite[Theorem~5.12]{PYZ}, their intersection is a finite intersection, so
$K'/\widetilde K'$ is finite. Then $K/\widetilde K$ is bounded and $\dim(\widetilde K)=\dim(K)$. Intersecting the growth chain with $\widetilde K$, we have a desired growth chain ending at $\widetilde K$.
\end{proof}

To obtain uniform results for definably amenable groups in $p$-adically closed fields, we study their analogues
over $N$, given by the first-order condition $\Good^*$. For the remainder of this subsection, fix a tuple $a$
from $N$ with $N\models\Good^*(a)$, and let $G=G_{N,a}$

\begin{lemma}\label{D:lem:semisimplebounded}
If $\operatorname{Lie}(G)$ is semisimple, then $G$ is bounded.
\end{lemma}

\begin{proof}
In $\mathcal U$-almost all factors, the corresponding
group $G'$ is definably amenable and
$\operatorname{Lie}(G')$ is semisimple.
Its adjoint homomorphism has finite kernel, since the differential $\operatorname{ad}$ is injective.
The Zariski closure of its image has semisimple
identity component.
Hence, the image is definably compact
\cite[Lemma~8.4]{PYZ}, and so is $G'$. Then $G'$ is bounded, and by {\L}o\'s's theorem, so is $G$.
\end{proof}

\begin{prop}\label{D:prop:completechains}
$G$ admits an $N$-definable complete growth chain with normal endpoint.
\end{prop}

\begin{proof}
We first prove the case when $G\leq\mathrm{GL}_m(N)$.
Let $\mathbf G=\operatorname{Zcl}(G)$, $\mathbf R$ be the solvable radical of $\mathbf G^\circ$. Let $W=G\cap\mathbf R(N)$.
The algebraic quotient map induces a definable homomorphism
$f:G\longrightarrow (\mathbf G/\mathbf R)(N)$
with kernel $W$. Since $G$ is open in $\mathbf G(N)$, the differential of $f$
at $1_G$ agrees with that of the algebraic quotient map, which is surjective by \cite[Chapter~II, Proposition~3.29]{Milne}. Thus,
$\operatorname{Lie}(f(G)) =\operatorname{Lie}(\mathbf G/\mathbf R)$ is semisimple. Since definable amenability is preserved under definable images in each $p$-adic factor, {\L}o\'s's theorem implies that
the parameters defining $f(G)$ in $N$ satisfy the $\operatorname{Grp}^*$-formula for the family of image groups. Then by Lemma~\ref{D:lem:semisimplebounded}, $f(G)$ is bounded. So we only need to deal with the part in $\mathbf R$.

Intersecting $G$ with the terms of the derived series of
$\mathbf R$, taken in reverse order, we have a series
\[
 \{1_G\}=V_0\lhd V_1\lhd\cdots\lhd V_s=W
\]
of subgroups normal in $G$, with each $V_{i+1}/V_i$
commutative. We construct a complete growth chain from the trivial growth chain. Suppose
that the current endpoint is $H\leq V_i$, where $H\lhd G$
and $V_i/H$ is bounded. If $V_{i+1}/H$ is bounded, we do not extend the growth chain at this step. Otherwise, we apply Lemma~\ref{D:lem:extension}, taking $V_{i+1}$
and $V_i$ for the groups denoted by $V$ and $W$ in the lemma, respectively. Its hypotheses hold since $V_i/H$ is bounded and $V_{i+1}/V_i$ is commutative. Since $V_{i+1}\lhd G$, the lemma allows us to extend the growth chain
to a new endpoint $\widetilde H\lhd G$ such that
$V_{i+1}/\widetilde H$ is bounded. Repeating this procedure along the series $(V_i)_i$, we obtain a growth chain ending at $K\lhd G$ with $W/K$ bounded. Since $G/W$ is
bounded, so is $G/K$, and hence the growth chain is now complete.

Now we prove the general case. 
Let $A=\ker\operatorname{Ad}$. In $\mathcal U$-almost all factors,
$Z(A')$ has finite index in $A'$
\cite[Lemmas~8.5 and~8.7]{ALH}, so $A/Z(A)$ is bounded. 
If $Z(A)$ is bounded, we let $H$ be trivial. 
If $Z(A)$ is unbounded, we apply Lemma~\ref{D:lem:extension} to the commutative normal subgroup $Z(A)$ of $G$.
In either case, we obtain a growth chain ending at $H\lhd G$ with $Z(A)/H$ bounded. Since $A/Z(A)$ and $Z(A)/H$ are bounded,  $A/H$ is bounded.

Now apply the linear construction to $\operatorname{Ad}(G)$,
whose images in almost all factors are definably amenable.
Pulling back the normal series used there along $\operatorname{Ad}$,
we have a series
\[
 A=V_0\lhd V_1\lhd\cdots\lhd V_s
\]
of normal subgroups of $G$ such that each $V_{i+1}/V_i$ is commutative and $G/V_s$ is bounded. Since $V_0/H$ is bounded, the induction
from the linear case works. It completes the proof.
\end{proof}

\subsection{Uniform Zariski closures}
\label{D:subsec:zariski}

We use the equality between definable dimension and the dimension of the Zariski closure, valid in $p$-adically closed fields and their ultraproducts. Together with compactness, we have a uniform degree
bound for equations defining Zariski closures and uniform codes for these closures.

\begin{lemma}\label{D:lem:zariski}
Fix an $L$-formula $\theta(z;w)$.
There is an integer $d$ such that, for every prime $p$,
every $M=M_p\models\pCF$, and every tuple $a$ from $M$, the Zariski closure of $\theta(M;a)$ is the common zero set of finitely many polynomials over $M$ of degree at most $d$.

Moreover, there are an integer $k$ and an $L$-formula
$\kappa(w;u)$ defining, in every such $M$, a function 
 $c_M:M^{|w|}\rightarrow M^k$
such that, for all $a,a'\in M^{|w|}$,
$$c_M(a)=c_M(a') \text{ iff } \operatorname{Zcl}(\theta(M;a)) = \operatorname{Zcl}(\theta(M;a')).$$
\end{lemma}

\begin{proof}
Suppose that no uniform degree bound exists.
For each $j\geq1$, pick $M_{p_j}$ and $a_j$ such that the Zariski closure of
 $X_j=\theta(M_{p_j};a_j)$
cannot be defined by polynomials of degree at most $j$.
After passing to a subsequence, assume that the primes are either constant or tend to infinity.
Let
 $N=\prod_j M_{p_j}/\mathcal U$, $
 a=(a_j)/\mathcal U$ and $X=\theta(N;a)$.
In either case, definable dimension in $N$ agrees with the dimension of the Zariski closure, as recalled above.

Let $Z_1,\ldots,Z_s$ be the geometric irreducible components of $\operatorname{Zcl}(X)$, all defined over $N$. Then
$X\subseteq\bigcup_{i=1}^s Z_i(N)$ and for every $i$, $\dim\bigl(X\cap Z_i(N)\bigr)=\dim Z_i$.

For each $i$, choose a finite set of equations over $N$ defining $Z_i$, and let $l$ be a common bound on their degrees.
For these fixed equation bounds, geometric irreducibility and the displayed conditions are first order. By {\L}o\'s's theorem, these properties hold for the corresponding algebraic sets in $\mathcal U$-almost all factors, so their union is exactly $\operatorname{Zcl}(X_j)$.
This union is defined by the products obtained by choosing one equation from each set.
These products have degree at most $sl$, contradicting the choice of $X_j$ for $j>sl$.
Thus, a uniform degree bound $d$ exists.

Fix an ordering of the monomials of degree at most $d$,
and identify polynomials with their coefficient vectors.
For each $a$, let
\[
 I_d(a)=
 \{P\in M[z]_{\leq d}:P(x)=0
       \text{ for every }x\in\theta(M;a)\}.
\]
These spaces form a uniformly $L$-definable family of vector subspaces of a coefficient space of fixed dimension. Clearly, the common zero set of $I_d(a)$ is precisely $\operatorname{Zcl}(\theta(M;a))$.

Using row reduction as in the proof of
Lemma~\ref{R:lem:linearalgebra}, let $c_M(a)$ be the tuple of entries of the unique reduced row-echelon matrix with row space $I_d(a)$, padded with zero rows to a fixed size.
Its graph is uniformly defined by an $L$-formula
$\kappa(w;u)$.

Then, for $a,a'\in M^{|w|}$, $\operatorname{Zcl}(\theta(M;a)) = \operatorname{Zcl}(\theta(M;a'))$ iff $I_d(a)=I_d(a')$ iff $c_M(a)=c_M(a')$.
\end{proof}

For subsets of $\mathrm{GL}_m(M)$, we apply the lemma
using the affine realization
\[
 g\longmapsto\bigl(g,\det(g)^{-1}\bigr).
\]
Thus, closures are taken inside the algebraic group $\mathrm{GL}_m$.
In particular, for a fixed definable family of subgroups of $\mathrm{GL}_m(M)$, the Zariski closures of these subgroups and their cosets admit uniform codes using only their defining parameters.

\section{Uniform dfg--fsg decompositions}
\label{sec:uniform-dfg-fsg}

In this section, we prove Theorem~\ref{main-H} and Theorem~\ref{main-HC}.

\subsection{Uniform normal dfg components}
\label{D:subsec:uniform-components}

We first show that complete growth chains can be constructed uniformly.

\begin{prop}\label{D:prop:finitechains}
There are finitely many $L$-formulas $\sigma_\ell(w;v)$ and $\psi_{\ell,i}(x;w,v)$ with $i\leq s_\ell$ such that for every $M=M_p\models\pCF$ and
every strongly admissible $a$ in $M$, there are $\ell$ and $b$
such that $M\models\sigma_\ell(a;b)$ and 
$(\psi_{\ell,i}(M;a,b))_{i\leq s_\ell}$ is a complete growth chain in $G_{M,a}$ with a normal endpoint.
\end{prop}

\begin{proof}
In each $M_p$ and tuple $a$ in $M_p$,  $\Good^*(a)$ means that $G_{M_p,a}$ is definably amenable. The dfg-fsg decomposition theorem \cite[Theorem~8.8]{PYZ}  and the structure theorem for dfg groups \cite[Theorem~5.6]{PYZ} give a complete growth chain with normal endpoint.
So there is a chain of formulas $(\psi_i)_{i\leq s}$ defining the complete growth chain. And there is $\sigma$ such that $\sigma(a)$ says that $(\psi_i)_{i\leq s}$ defines a desired growth chain for $G_{{M_p},a}$. Replacing the extra parameters  by a tuple of free variables $v$, we obtain $L$-formulas $\sigma(w;v)$ and $\psi_i(x;w,v)$. We enumerate all these possible sequences of $L$-formulas as $(\sigma_\ell, \psi_{\ell,i}(x;w,v_\ell):0\leq i\leq s_\ell)$.

Now it suffices to show that there is $k$ such that $$M_p\models \forall w(\Good^*(w)\rightarrow \bigvee_{\ell\leq k}\exists v_\ell \sigma_\ell(w,v_\ell))$$
for every prime $p$.
If not, for each $j$, there are prime $p_j$ and $a_j\in M_{p_j}$ such that $M_{p_j}\models \Good^*(a_j)\wedge \bigwedge_{\ell\leq j}\neg\exists v_\ell\sigma_\ell(a_j;v_\ell)$. Let $N=\prod_j M_{p_j}/\mathcal U$ be a nonprincipal ultraproduct and let
$a=(a_j)/\mathcal U$. Passing to a subsequence, we may assume that either $p_j$ is
constant or $p_j\to\infty$.
In either case, $G_{N,a}$ admits a complete growth chain with normal endpoint (by Proposition~\ref{D:prop:completechains} when $p_j\rightarrow\infty$). So there is $k$ such that $N\models \Good^*(a)\wedge \exists v_k \sigma_k(a,v_k)$. Then there is $j>k$ such that $M_{p_j}\models \Good^*(a_j)\wedge \exists v_k \sigma_k(a_j,v_k)$, which is a contradiction.
\end{proof}

Now we prove Theorem~\ref{main-H}.

\begin{proof}[Proof of Theorem~\ref{main-H}]

We follow the notations in Proposition~\ref{D:prop:finitechains} and assume $\ell\leq k$. Then for every $M_p$ and tuples $a$ and $b$ in $M_p$,  if $\sigma_\ell(a,b)$ holds, then $\psi_{\ell,s_\ell}(M_p;a,b)$ is a normal dfg component of $G_{M_p,a}$. To remove the extra parameters $b$, we collect all these normal dfg components and let
$$\mathcal Y_a:=\bigcup_{\ell\leq k}\{\psi_{\ell,s_\ell}(M_p;a,b): \sigma_\ell(a,b) \text{ holds}\}$$
Then $\mathcal Y_a$ is a definable family of normal dfg components of $G_{M_p,a}$.
By \cite[Theorem~5.12]{PYZ}, the intersection $\bigcap \mathcal Y_a$
is a finite intersection of members of $\mathcal Y_a$.
Repeated application of \cite[Corollary~8.9]{PYZ} shows that
this finite intersection is again a dfg component.
Hence, $\bigcap\mathcal Y_a$ is a normal dfg component of $G_{M_p,a}$.

Now we let
$$\delta_H(x;w):= \bigwedge_{\ell\leq k}\forall v
 \bigl(\sigma_\ell(w;v)\rightarrow
       \psi_{\ell,s_\ell}(x;w,v)\bigr).$$
Then $\delta_H(x;w)$ is an $L$-formula. It is easy to check that $\delta_H(x;a)$  defines $\bigcap \mathcal Y_a$ in every $M_p$ and strongly admissible $a$.
\end{proof}

\subsection{Uniform definable quotient maps}
\label{D:subsec:quotients}

We first prove the linear case of Theorem~\ref{main-HC}.
In this case, the codes for Zariski closures from
Subsection~\ref{D:subsec:zariski} give the quotient map directly.

\begin{proof}[Proof of the linear case of Theorem~\ref{main-HC}]
Suppose that the given family consists of subgroups of $\mathrm{GL}_m$ for a fixed $m$. For $M=M_p\models\pCF$ and
a strongly admissible tuple $a$ from $M$, let
 $H_{M,a}=\delta_H(M;a)$.
By Lemma~\ref{D:lem:zariski}, there is an $L$-formula such that $\rho_{M}$ gives the code of $\Zcl(gH_{M,a})$ for $g\in G_{M,a}$. We let $\rho_{M,a}(g):=\rho_M(g,a)$.
Since the left translation by $g$
is an algebraic map of $\mathrm{GL}_m$,
$\Zcl(gH_{M,a})=g\Zcl(H_{M,a})$. Hence, for $g,g'\in G_{M,a}$,
$\rho_{M,a}(g)=\rho_{M,a}(g')$ iff $g\Zcl(H_{M,a})=g'\Zcl(H_{M,a})$ iff $g^{-1}g'\in G_{M,a}\cap \Zcl(H_{M,a})(M)$.

Let $K_{M,a}= G_{M,a}\cap \Zcl(H_{M,a})(M)$. Since $H_{M,a}$ is normal in $G_{M,a}$, so is $K_{M,a}$. The fibers of
$\rho_{M,a}$ are precisely the left cosets of $K_{M,a}$.
Since $\dim H_{M,a}=\dim\Zcl(H_{M,a})$, the subgroup $H_{M,a}$ is open in $K_{M,a}$.
As $H_{M,a}$ has dfg, it has finite index in $K_{M,a}$
\cite[Corollary~4.4]{AGJ}. Thus, $K_{M,a}$ is a normal dfg component. Let $C_{M,a}:=\rho_{M,a}(G_{M,a})$ and define the multiplication as
\[
 \rho_{M,a}(g)\rho_{M,a}(g')=\rho_{M,a}(gg').
\]
This multiplication is well defined and $C_{M,a}$ has fsg. Then we have dfg-fsg decompositions uniformly
$$1\longrightarrow K_{M,a}\longrightarrow G_{M,a}\stackrel{\rho_{M,a}}{\longrightarrow} C_{M,a}\longrightarrow 1$$
as desired without extra parameters. 
\end{proof}

For a general group presentation, we construct the quotient along a growth chain. At each step, a one-dimensional orbit gives rise to a curve in field coordinates. Coding the curve
and its limits gives a map on the coset space with finite fibers. We then use the differential of this map on $G$ to
obtain a group quotient whose kernel is a finite extension of the endpoint.

For the next three lemmas, let $N=\prod_j M_{p_j}/\mathcal U$
with $p_j\to\infty$, as in Subsection~\ref{D:subsec:auxiliary}.
Work in a sufficiently saturated elementary extension $N'\succeq N$.

\begin{lemma}\label{D:lem:curve-types}
Let $C$ be an algebraic curve over $N'$ and let
$c\in C(N')$. There are at most $2^{\aleph_0}$ complete $L$-types over $N'$ concentrating on $C\setminus\{c\}$
and infinitesimally close to $c$.
\end{lemma}

\begin{proof}
We first bound the number of complete $L$-types over $N'$ of nonzero infinitesimals at $0$. Let $u$ realize such a type, so $v(u)>\gamma$ for every $\gamma\in\Gamma_{N'}$. For every nonzero polynomial $P$ over $N'$, the leading term of $P(u)$ is that of its lowest
nonzero monomial. By \cite[Propositions~2.2.2 and~2.3.2]{AcostaHils},
the type of $u$ is determined by its power cosets and the congruences of $v(u)$. For each $n$, there are finitely many
$n$-th-power cosets, since the residue field is pseudofinite of characteristic zero, the principal units are divisible,
and $\Gamma_{N'}/n\Gamma_{N'}$ is finite. This gives at most $2^{\aleph_0}$ types.

By Noether normalization, there is a finite-to-one regular map $t:U\to\mathbb A^1$ over $N'$ on an affine neighbourhood $U$
of $c$. We may assume that $t(c)=0$.
Since $t^{-1}(0)$ is finite, every $x\neq c$ infinitesimally close to $c$ maps to a nonzero infinitesimal. Each complete type of $t(x)$ over $N'$ has only finitely many lifts, so
the preceding bound applies.
\end{proof}

\begin{lemma}\label{D:lem:chain-cosets}
Let $G$ be a definable group and let $H$ be the endpoint of
a growth chain in $G$. Then there are a definable map
$f:G\to(N')^r$ and a natural number $d$ such that $f$ is
constant on each left $H$-coset and each fiber contains at
most $d$ such cosets.
\end{lemma}

\begin{proof}
We write the growth chain as
\[
 \{1_G\}=H_0\lhd H_1\lhd\cdots\lhd H_s=H.
\]
We construct the maps $f_i$ by induction on $i$.
Let $f_0$ be the identity map on $G$. Suppose that $f_i:G\to(N')^{r_i}$ is constant on left $H_i$-cosets, with at most $d_i$ such cosets in each fiber. Since $H_i\lhd H_{i+1}$,
the group $H_{i+1}/H_i$ acts on $G/H_i$ on the right by
$(gH_i)\cdot(hH_i)=ghH_i$. Its orbits are the fibers of
$G/H_i\to G/H_{i+1}$.

For $g\in G$, let $Y_g=f_i(gH_{i+1})$. Each fiber of $f_i$ is a union of at most $d_i$ left $H_i$-cosets. Applying the fiber dimension theorem to the restriction of $f_i$ to
$gH_{i+1}$, we obtain
\[
 \dim Y_g=\dim H_{i+1}-\dim H_i=1.
\]

Let $(E_\lambda)_\lambda$ be a
$\Gamma$-exhaustion of $G$. Let 
$$ Y_{g,\lambda}
 =f_i\bigl(gH_{i+1}\setminus g(E_\lambda\cap H_{i+1})H_i\bigr)$$
and $Z_g:=\bigcap_{\lambda} \cl(Y_{g,\lambda})$
where the closure is taken in the valuation topology in the projective space
$\mathbb P^{r_i}(N')$. Since translates of bounded sets are bounded, $Z_g$ is independent
of the choice of representative $g$ of the coset $gH_{i+1}$.

Since $H_{i+1}/H_i$ is unbounded, the sets $Y_{g,\lambda}$ are
nonempty and decrease with $\lambda$. By definable compactness of $\mathbb P^{r_i}(M_{p_j})$, the
corresponding closures have nonempty intersection in
$\mathcal U$-almost all factors. By {\L}o\'s's theorem and
$N\preceq N'$, it follows that $Z_g\neq\emptyset$.

\begin{claim}
The set $Z_g$ is finite.
\end{claim}

\begin{claimproof}
Each fiber of $f_i$ in $gH_{i+1}$ consists of at most $d_i$
left $H_i$-cosets, all contained in $g(E_\lambda\cap H_{i+1})H_i$
for sufficiently large $\lambda$. Hence,
$\bigcap_\lambda Y_{g,\lambda}=\emptyset$. 

Suppose that $Z_g$ is infinite. Then there are pairwise distinct $c_n\in Z_g$ for $n\in\mathbb N$. Since the sets $Y_{g,\lambda}$
decrease and have empty intersection, the saturation of $N'$ gives
a common $\lambda$ such that $c_n\in \cl(Y_{g,\lambda})\backslash Y_{g,\lambda}$ for every $n$. 
By \cite[Proposition~2.18]{ALH}, $\dim (\cl(Y_{g,\lambda})\backslash Y_{g,\lambda}) < \dim Y_{g,\lambda}\leq 1$. Hence, $\cl(Y_{g,\lambda})\backslash Y_{g,\lambda}$ is finite.
This is a contradiction.
\end{claimproof}

Let $f_{i+1}(g)$ be the pair consisting of the code for  the affine Zariski closure $\Zcl(Y_g)$ from Lemma~\ref{D:lem:zariski} and the code for the finite set $Z_g$ given by symmetric polynomials.
Then $f_{i+1}$ is definable and constant on each left $H_{i+1}$-coset.

Suppose that $f_{i+1}(g)=f_{i+1}(g')$. Then
$\Zcl(Y_g)=\Zcl(Y_{g'})$ and $Z_g=Z_{g'}$. Let $C$ be the Zariski closure of $\Zcl(Y_g)$ in $\mathbb P^{r_i}$ and fix some $c\in Z_g$.
For every such $g'$, we then have $Y_{g'}\subseteq C(N')$ and
$c\in Z_{g'}$.

By Lemma~\ref{D:lem:curve-types}, there are at most $2^{\aleph_0}$ complete types on $C\setminus\{c\}$ over $N'$ infinitesimally
close to $c$. Since each fiber of $f_i$ is a union of at most $d_i$ left $H_i$-cosets, each such type can contain $Y_{g'}$ for at most
$d_i$ distinct cosets $g'H_{i+1}$. Thus, the fiber contains at most $d_i\,2^{\aleph_0}$ cosets. These cosets form a definable subset of $G/H_{i+1}$, so the saturation of $N'$ implies that there are only finitely many. By saturation again, there is a common finite bound $d_{i+1}$ for all fibers. This completes the induction.
\end{proof}

\begin{lemma}\label{D:lem:real-group-quotient}
Let $G$ be a definable group, and let $H$ be
a definable normal subgroup. Suppose that there is a definable map
$f:G\rightarrow (N')^r$ such that each fiber is a union of at most $d$
left $H$-cosets, for some $d$.
Then there are a definable group $C$ and a definable surjective
homomorphism $\rho:G\rightarrow C$ such that
$\ker \rho$ is a finite extension of $H$.
\end{lemma}

\begin{proof}
Let $Q=G/H$ and let $\pi:G\to Q$ be the quotient map.  At each point $g$ where $f$ is strictly differentiable, let $J_f(g)$ be the linear map
$\mathrm{d}(f\circ R_g)_{1_G}: \Lie(G)\rightarrow (N')^r$ where $R_g: G\rightarrow G: x\mapsto xg$.

By \cite[Proposition~5.7]{ALH}, $f$ is strictly differentiable on a definable open dense subset of $G$. Its differential is continuous there, so its rank is locally constant on a definable
open dense subset $U$ of $G$. Since each fiber of $f$ is a finite union of left $H$-cosets, every nonempty relatively open definable
subset of a fiber has dimension $\dim H$. By \cite[Proposition~4.11]{ALH}, we have that, for $g\in U$,
$$\operatorname{rank}J_f(g)=\dim G-\dim H.$$
Since $H$ is normal in $G$, the map $f\circ R_g$ is constant on $H$, so $\Lie(H)\subseteq\ker J_f(g)$. It follows
that $\ker J_f(g)=\Lie(H)$ for every $g\in U$.
Since $f\circ R_h=f$ for every $h\in H$, we may replace $U$ by $UH$ and assume that $U$ is right $H$-invariant.

Since $H$ is normal in $G$, the adjoint action $\operatorname{Ad}(G)$ preserves $\Lie(H)$. Let $\alpha$ denote the induced action
of $G$ on $\Lie(G)/\Lie(H)$.
For $h\in H$, the map $x\mapsto hxh^{-1}x^{-1}$ takes values
in $H$, so its differential $\operatorname{Ad}(h)-\mathrm{id}$
at $1_G$ has image in $\Lie(H)$. Hence, $H$ acts trivially on $\Lie(G)/\Lie(H)$.

Then for $g\in U$, let
\[
 F_0(g)=\bigl(f(g),J_f(g),\alpha(g)\bigr).
\]
where $\alpha(g)$ is written in a fixed
basis. The map $F_0$ is right $H$-invariant.  Let
\[
 A_0=\{ag^{-1}g'a^{-1}:a\in G \text{ and } g,g'\in U \text{ with } F_0(g)=F_0(g')\}.
\]
and let $A=\pi(A_0)$. Since $A_0$ is preserved under conjugation in $G$, its image $A$ is preserved under conjugation in $Q$.

\begin{claim}
$A$ is finite.
\end{claim}

\begin{claimproof}
Let $Z=\{(g,g')\in U\times U:F_0(g)=F_0(g')\}$ and $\Phi(a,g,g')=ag^{-1}g'a^{-1}$. Clearly, $A_0=\Phi(G\times Z)$.
Working in charts, partition $Z$ into finitely many definable
\'etale pieces by \cite[Propositions~2.15 and~5.9]{ALH}.
For each piece $S$, the inclusion $\ell:S\hookrightarrow U^2$ is
strictly differentiable on a definable open dense subset
$V$ by \cite[Proposition~5.7]{ALH}.

At $(g,g')\in V$, let $(\xi,\xi')\in\Lie(G)\times\Lie(G)$
represent an arbitrary vector in the image of
$\mathrm{d}\ell_{(g,g')}$ in right translation coordinates. 
Differentiating $f(g)=f(g')$ along
$V$ gives $J_f(g)\xi=J_f(g')\xi'$. Since $J_f(g)=J_f(g')$,
we obtain $\xi'-\xi\in\Lie(H)$. Let $z=\Phi(a,g,g')=ag^{-1}g'a^{-1}$ and $w=g^{-1}g'$. Then $z=awa^{-1}$.
In right translation coordinates, the differential of multiplication
at $(u,v)$ sends $(X,Y)$ to $X+\operatorname{Ad}(u)Y$, and that
of inversion at $u$ sends $X$ to $-\operatorname{Ad}(u^{-1})X$.
Thus, the tangent vector at $w$ is
$\omega=\operatorname{Ad}(g^{-1})(\xi'-\xi)$. Writing $\zeta$
for the tangent vector at $a$, we obtain at $z$,
\[
\begin{aligned}
 \zeta+\operatorname{Ad}(a)\omega
 -\operatorname{Ad}(aw)\operatorname{Ad}(a^{-1})\zeta
 &=(\mathrm{id}-\operatorname{Ad}(z))\zeta
   +\operatorname{Ad}(ag^{-1})(\xi'-\xi).
\end{aligned}
\]
Since $\alpha(g)=\alpha(g')$, we have
$\alpha(z)=\alpha(a)\alpha(g)^{-1}\alpha(g')\alpha(a)^{-1}=\mathrm{id}$.
Thus, $\operatorname{Ad}(z)$ acts trivially on $\Lie(G)/\Lie(H)$,
so $(\mathrm{id}-\operatorname{Ad}(z))\zeta\in\Lie(H)$.
Moreover, $\xi'-\xi\in\Lie(H)$, and normality of $H$ implies
that $\operatorname{Ad}(ag^{-1})$ preserves $\Lie(H)$. Hence
$\operatorname{Ad}(ag^{-1})(\xi'-\xi)\in\Lie(H)$ as well.
Thus, the differential of $\Phi$ on $G\times V$ has rank at most $\dim H$, so $\dim\Phi(G\times V)\leq\dim H$ by
\cite[Proposition~5.21]{ALH}.

By continuity, $
 \Phi(G\times S)\subseteq\cl(\Phi(G\times V))$.
Taking closures preserves dimension \cite[Proposition~2.18]{ALH}.
Since the partition is finite, $\dim A_0\leq\dim H$.
Right $H$-invariance of $U$ and $F_0$, together with normality of $H$, shows that $A_0$ is a union of left $H$-cosets. Hence, every
nonempty fiber of $f|_{A_0}$ has dimension $\dim H$. The fiber dimension theorem gives $\dim f(A_0)=0$, so $f(A_0)$ is finite.
As each fiber contains at most $d$ left $H$-cosets, $A=\pi(A_0)$ is finite.
\end{claimproof}

By \cite[Lemma~6.2]{ALH}, finitely many left translates $a_1U,\ldots,a_kU$ cover $G$. For $g\in G$,
take the least $j$ with $g\in a_jU$ and let
\[
 F_1(\pi(g))=\bigl(j,F_0(a_j^{-1}g)\bigr).
\]
This map is well defined and its fibers have size at most $d$.
If $F_1(\pi(g))=F_1(\pi(g'))$, the same index $j$ occurs for both,
so $\pi(g^{-1}g')\in A$.

Let $P$ be the subgroup of $Q$ generated by $A$. The conjugation action of $P$ on the
finite set $A$ has kernel $Z(P)$.
Then $P/Z(P)$ is finite, and by Schur's theorem, the derived subgroup $[P,P]$ is finite. Let $K$ be the inverse image in $P$ of the torsion subgroup of $P/[P,P]$. Since $K/[P,P]$ is the torsion subgroup of the
finitely generated abelian group $P/[P,P]$, it is finite. Since $[P,P]$ is finite, $K$ is finite. Moreover, since $K$ is characteristic in $P$ and $P$ is normal in $Q$,
$K$ is a finite normal subgroup of $Q$.

Let $W=\{u^j:u\in A \text{ and }0\leq j<d\}$. Fix an order of $W$ and define a map on $Q$ as
\[
 F_2(q)=\bigl(F_1(qw)\bigr)_{w\in W}.
\]
Suppose $F_2(\pi(g))=F_2(\pi(g'))$, and let
$u=\pi(g^{-1}g')$. Since $1_Q\in W$, equality of $F_2(\pi(g))$ and $F_2(\pi(g'))$
implies $F_1(\pi(g))=F_1(\pi(g'))$, so $u\in A$. Since
$\pi(g')=\pi(g)u$, comparing the coordinates indexed by
$u^j$ for $0\leq j<d$, we have
\[
 F_1(\pi(g))=F_1(\pi(g)u)=\cdots=F_1(\pi(g)u^d).
\]
Two of these $d+1$ arguments coincide by the fiber bound for $F_1$. Hence, $u$ has order at most $d$ and therefore belongs to $K$.

Finally, we define a map $\rho$ on $G$ as
\[
 \rho(g)=\text{ the code of the finite set }\{F_2(\pi(g)k): k\in K\}.
\]
Clearly, $\rho(g)$ depends only on $\pi(g)K$. Conversely, if
$\rho(g)=\rho(g')$, then $F_2(\pi(g))=F_2(\pi(g')k)$ for some
$k\in K$. By the preceding paragraph, $\pi(g)K=\pi(g')K$.
Thus multiplication on $Q/K$ induces a definable group law
on $C=\rho(G)$, making $\rho:G\to C$ a surjective homomorphism
with kernel $\pi^{-1}(K)$, a finite extension of $H$.
\end{proof}

Now we prove the general case of Theorem~\ref{main-HC}.

\begin{proof}[Proof of the general case of Theorem~\ref{main-HC}]
In each $M_p$ and strongly admissible $a$ in $M_p$, \cite[Theorem~8.8]{PYZ} gives a dfg-fsg decomposition of $G_{M_p,a}$. Moreover, we can assume that the dfg part
contains $H_{M_p,a}:=\delta_H(M_p;a)$ as a subgroup of finite index, where $\delta_H(x;w)$ is the formula from Theorem~\ref{main-H}.
Thus, there are an $L$-formula $\chi(x,y;w,v)$, a natural number $d$, and a tuple $b$ such that $\chi(x,y;a,b)$ defines the graph of the quotient map $\rho$ of a dfg-fsg decomposition of $G_{M,a}$  with $H_{M_p,a}\leq\ker\rho$ and
$|\ker\rho/H_{M_p,a}|\leq d$.

Since the dfg condition on the kernel and the fsg condition on the image are definable in families by \cite[Theorem~5.9]{PYZ} and \cite[Theorem~1.2 and Proposition~4.1]{JohnsonFsg}, respectively, there is an $L$-formula $\sigma(w,v)$ and $\sigma(a,b)$ expresses the first-order conditions ensuring that $\chi(x,y;a,b)$ defines the graph of the quotient map in a dfg-fsg decomposition such that the subgroup defined by $\delta_H$ has index at most $d$ in the kernel of this map.

We enumerate all these possible pairs of $L$-formulas and the corresponding bounds as $(\chi_\ell(x,y_\ell;w,v_\ell), \sigma_\ell(w,v_\ell),d_\ell)$.

To obtain dfg-fsg decompositions uniformly in which the group defined by $\delta_H$ has uniformly bounded index in the kernel, it suffices to show that there is $k$ such that $$M_p\models \forall w(\Good^*(w)\rightarrow \bigvee_{\ell\leq k}\exists v_\ell \sigma_\ell(w,v_\ell))$$
for every prime $p$.
If not, for each $j$, there are prime $p_j$ and $a_j\in M_{p_j}$ such that $M_{p_j}\models \Good^*(a_j)\wedge \bigwedge_{\ell\leq j}\neg\exists v_\ell\,\sigma_\ell(a_j;v_\ell)$. Let $N=\prod_j M_{p_j}/\mathcal U$ be a nonprincipal ultraproduct and let
$a=(a_j)/\mathcal U$. Passing to a subsequence, we may assume that either $p_j$ is
constant or $p_j\to\infty$.

If $p_j=p$ for all $j$, then $N\models\pCF$. Since $a$ is
strongly admissible, the fixed-prime case gives an index $\ell$
with $d_\ell=1$ such that $N\models\exists v_\ell\,\sigma_\ell(a,v_\ell)$.

Suppose that $p_j\to\infty$, and let $N'\succeq N$ be sufficiently
saturated. In each factor, $\delta_H(M_{p_j};a_j)$ has finite
index in the endpoint of a chain from
Proposition~\ref{D:prop:finitechains}. Intersecting the chain
with this subgroup and applying {\L}o\'s's theorem and
elementarity gives a complete growth chain ending at
$\delta_H(N';a)$. By Lemmas~\ref{D:lem:chain-cosets}
and~\ref{D:lem:real-group-quotient}, there is a definable
quotient map $\rho$ on $G_{N',a}$ whose kernel contains
$\delta_H(N';a)$ with finite index and whose quotient is bounded.
Taking $d=|\ker\rho/\delta_H(N';a)|$, we obtain an index $\ell$
such that $N'\models\exists v_\ell\,\sigma_\ell(a,v_\ell)$.
By elementarity, $N\models\exists v_\ell\,\sigma_\ell(a,v_\ell)$.

In either case, {\L}o\'s's theorem gives
$M_{p_j}\models\exists v_\ell\,\sigma_\ell(a_j,v_\ell)$
for $\mathcal U$-almost all $j$, contradicting the choice of
$a_j$ for $j\geq\ell$.

It remains to show that the kernel is independent of $b$. For each
strongly admissible $a$, let $K_{M_p,a}$ be the intersection
of the kernels of all maps given by $\chi_\ell(x,y_\ell;a,b)$
with $\ell\leq k$ and $M_p\models\sigma_\ell(a,b)$.  It is
uniformly definable over $a$ and is the intersection of at most
$d=\max_{\ell\leq k}d_\ell$ of these kernels, since each strict
intersection decreases the index over $H_{M_p,a}$. The induced map into the Cartesian product of the corresponding quotient groups has kernel $K_{M_p,a}$. Taking
its image as the target gives the required dfg-fsg decomposition.
\end{proof}

\section*{Declaration of AI use}

The author used ChatGPT (OpenAI) for assistance with language editing and for discussing and checking mathematical arguments during the preparation of this manuscript. All mathematical statements, arguments, and proofs included in the final version were independently verified by the author, who takes full responsibility for the content of the article.

\end{document}